\documentclass[11pt, twoside]{article}
\usepackage{latexsym}
\usepackage{amsmath}
\usepackage{amssymb}
\usepackage[all]{xy}
\usepackage{amsfonts}
\usepackage{verbatim}
\usepackage{amsthm}
\usepackage{bm}
\usepackage{mathrsfs}
\usepackage{epsfig}
\usepackage{xy}
\usepackage{array}
\usepackage{stmaryrd}
\usepackage{graphicx,color}
\usepackage{xcolor}
\usepackage{tikz}
\usepackage{xcolor}
\usepackage[colorlinks=true,linkcolor=blue,citecolor=blue]{hyperref}
\usetikzlibrary{arrows,calc}
\usepackage{etex}
\usepackage{mathdots}
\usepackage{float}
\usepackage{graphics}
\usepackage{pdflscape}
\usepackage{CJK}
\usepackage{anysize,hyperref}
\input xypic
\xyoption{all}
\usepackage{extarrows}
\usepackage[perpage,symbol]{footmisc}
\usepackage{setspace}
\usepackage{enumitem}
\setenumerate[1]{itemsep=0pt,partopsep=0pt,parsep=\parskip,topsep=5pt}
\setitemize[1]{itemsep=0pt,partopsep=0pt,parsep=\parskip,topsep=5pt}
\setdescription{itemsep=0pt,partopsep=0pt,parsep=\parskip,topsep=5pt}

\def\A{\mathcal{A}}

\def\T{\mathcal{T}}

\def\C{\mathscr{C}}
\def\F{\mathbb{S}}

\def\dr{\ar@{->}[r]}

\def\X{\mathscr{X}}

\def\add{\mbox{add}}

\def\Hom{\mbox{Hom}}

\def\rad{\mbox{rad}}
\def\End{\mbox{End}}
\def\Ker{\mbox{Ker}}

\renewcommand{\diagram}[3]{\matrix (#1) [matrix of math nodes,row
  sep={#2},column sep={#3},text height=1.5ex,text depth=0.25ex]}
\usepackage{pgf,tikz}

  \usetikzlibrary{arrows}
  \usetikzlibrary{decorations.markings}
  \usetikzlibrary{shapes}
  \usetikzlibrary{matrix}
  \usepackage{wasysym}
\begin{document}
\baselineskip=15pt
\title{\Large{\bf Higher-dimensional Auslander-Reiten theory on \\[2mm] $\bm{(d+2)}$-angulated categories$^\bigstar$\footnotetext{\hspace{-1em}$^\bigstar$This work was supported by the Scientific Research Fund of Hunan Provincial Education Department (Grant No. 19B239) and the National Natural Science Foundation of China (Grant No. 11901190).}}}
\medskip
\author{Panyue Zhou}

\date{}

\maketitle
\def\blue{\color{blue}}
\def\red{\color{red}}

\newtheorem{theorem}{Theorem}[section]
\newtheorem{lemma}[theorem]{Lemma}
\newtheorem{corollary}[theorem]{Corollary}
\newtheorem{proposition}[theorem]{Proposition}
\newtheorem{conjecture}{Conjecture}
\theoremstyle{definition}
\newtheorem{definition}[theorem]{Definition}
\newtheorem{question}[theorem]{Question}
\newtheorem{remark}[theorem]{Remark}
\newtheorem{remark*}[]{Remark}
\newtheorem{example}[theorem]{Example}
\newtheorem{example*}[]{Example}
\newtheorem{condition}[theorem]{Condition}
\newtheorem{condition*}[]{Condition}
\newtheorem{construction}[theorem]{Construction}
\newtheorem{construction*}[]{Construction}

\newtheorem{assumption}[theorem]{Assumption}
\newtheorem{assumption*}[]{Assumption}

\baselineskip=17pt
\parindent=0.5cm

\begin{abstract}
\baselineskip=16pt
Let $\C$ be a $(d+2)$-angulated category with $d$-suspension
functor $\Sigma^d$. Our main results show that every Serre functor on $\C$ is a $(d+2)$-angulated functor.
We also show that $\C$ has a Serre functor $\mathbb{S}$ if and only if $\C$ has Auslander--Reiten $(d+2)$-angles.
Moreover, $\tau_d=\mathbb{S}\Sigma^{-d}$ where $\tau_d$ is $d$-Auslander-Reiten translation.
These results generalize work by Bondal-Kapranov and  Reiten-Van den Bergh.
As an application, we prove that for a strongly functorially finite subcategory $\X$ of $\C$, the quotient category
$\C/\X$ is a $(d+2)$-angulated category if and only if $(\C,\C)$ is an $\X$-mutation pair, and if and
only if $\tau_d\X =\X$.\\[0.5cm]
\textbf{Keywords:} $(d+2)$-angulated categories; Auslander-Reiten $(d+2)$-angles; Serre functor\\[0.2cm]
\textbf{ 2020 Mathematics Subject Classification:} 16G70; 18G80
\medskip
\end{abstract}

\pagestyle{myheadings}
\markboth{\rightline {\scriptsize   Panyue Zhou}}
         {\leftline{\scriptsize  Higher-dimensional Auslander-Reiten theory on $(d+2)$-angulated categories}}

\section{Introduction}
Unless otherwise specified $k$ will be an algebraically closed field and
all categories in this article will be $k$-linear Hom-finite Krull-Schmidt.
Let $\C$ be an additive category. We say that $\C$ has a Serre functor $\mathbb{S}$, that is, an auto-equivalence for which there exists a natural equivalence
$$\Hom_{\C}(X,Y)\simeq\Hom_{\C}(Y, \mathbb{S}X)^{\ast}$$
for any $X,Y\in\C$, where $(-)^{\ast}=\Hom_k (-, k)$ is the $k$-linear duality functor.

Bondal and Kapranov showed the following result. Van den Bergh has also given other proof of this result in the appendix of this article \cite[Theorem A.4.4]{b}.
\begin{theorem}\emph{\cite[Proposition 3.3]{bk}}
Let $\C$ be a triangulated category. Then every Serre functor is a
triangulated functor.
\end{theorem}

Reiten and Van den Bergh gave a connection between Auslander-Reiten triangles and Serre functor.
\begin{theorem}\emph{\cite[Theorem I.2.4]{rv}}
Let $\C$ be a triangulated category with shift functor $\Sigma$. Then
$\C$ has a Serre functor $\mathbb{S}$ if and only if $\C$ has Auslander-Reiten triangles.
Moreover, $\tau=\mathbb{S}\Sigma^{-1}$ where $\tau$ is  the Auslander-Reiten translation.
\end{theorem}

Quotient categories come in a number of different flavours. The one to be considered here is probably the most basic. Let $\C$ be an additive category with a subcategory $\X$. For two objects $A,B\in\C$, denote by $\X(A,B)$ all the morphisms from
$A$ to $B$ which factor through an object of $\X$. Then the quotient category $\C/\X$ has
the same objects as $\C$, and its morphism spaces are defined by
$$\Hom_{\C/\X}(A,B)=\Hom_{\C}(A,B)/\X(A,B).$$

J{\o}rgensen gave a necessary and sufficient condition for quotient categories to be triangulated categories.
\begin{theorem}\emph{\cite[Theorem 3.3]{j}}
Let $\C$ be a triangulated category with Serre functor and $\X$ be a functorially finite subcategory of $\C$. Then
the quotient category $\C/\X$ is a triangulated category if and only if $\tau\X=\X$.
\end{theorem}

In \cite{gko}, Geiss, Keller and Oppermann introduced $(d+2)$-angulated categories. These are
generalizations of triangulated categories, in the sense that triangles are replaced by
$(d+2)$-angles, that is, morphism sequences of length $d+2$. In particular, a $3$-angulated category is precisely
a triangulated category.
The notion of Auslander-Reiten $(d+2)$-angles in a special $(d+2)$-angulated category was first introduced by
Iyama and Yoshino \cite{iy}, and then was generalized by Fedele \cite{f} to arbitrary $(d+2)$-angulated category.

In this article we will generalize these results into the higher homological case.
Moreover, our proof is not far from the usual triangulated case. We
hope that our work would motivate further study on $(d+2)$-angulated categories.

Our first main result is the following.
\begin{theorem}\emph{(see Theorem \ref{main1} for details)}
Let $\C$ be a $(d+2)$-angulated category with  $d$-suspension
functor $\Sigma^d$ and right Serre functor  $\F$.
Then there exists a  natural isomorphism $\zeta\colon \F\Sigma^d\to \Sigma^d\F$ such that
$\F\colon \C\to\C$ is a $(d+2)$-angulated functor.
\end{theorem}

Using the first main result, we get our second main result which is as follows.

\begin{theorem}\emph{(see Theorem \ref{main2} for details)}
Let $\C$ be a $(d+2)$-angulated category with a $d$-suspension
functor $\Sigma^d$.
Then $\C$ has Auslander-Reiten $(d+2)$-angles if and only if
$\C$ has a Serre functor $\mathbb{S}$. Moreover $\tau_d=\mathbb{S}\Sigma^{-d}$ where $\tau_d$ is
the $d$-Auslander-Reiten translation.
\end{theorem}

As an application of our second main result, we have the following third main result.

\begin{theorem}\emph{(see Theorem \ref{main3} for details)}
Let $\C$ be a $(d+2)$-angulated category with $d$-suspension
functor $\Sigma^d$ and a Serre functor $\mathbb{S}$, and $\X$ be a strongly functorially finite
subcategory of $\C$. Then the following statements are equivalent:
\begin{itemize}
\item[\rm (1)] $(\C,\C)$ is an $\X$-mutation pair.

\item[\rm (2)] The quotient category $\C/\X$ is a $(d+2)$-angulated category.

\item[\rm (3)] $\tau_d\X=\X$.
\end{itemize}
\end{theorem}

The article is organized as follows: In Section 2, we review some elementary definitions
that we need to use, including  $(d+2)$-angulated categories and Auslander-Reiten $(d+2)$ angles.
In Section 3,  we show the first main result. In Section 4, we prove the second main result and give an example illustrating it. In Section 5, we show the third main result and give an example illustrating it.

\section{Preliminaries}
\subsection{Right $(d+2)$-angulated categories}

Let $\C$ be an additive category with an endofunctor $\Sigma^d\colon\C\rightarrow\C$ ($\Sigma^d$ is called the $d$-suspension functor), and $d$ be an integer greater than or equal to one. A $(d+2)$-$\Sigma^d$-$sequence$ in $\C$ is a sequence of morphisms
$$A_0\xrightarrow{f_0}A_1\xrightarrow{f_1}A_2\xrightarrow{f_2}\cdots\xrightarrow{f_{d-1}}A_d\xrightarrow{f_d}A_{d+1}\xrightarrow{f_{d+1}}\Sigma^d A_0.$$
Its {\em left rotation} is the $(d+2)$-$\Sigma^d$-sequence
$$A_1\xrightarrow{f_1}A_2\xrightarrow{f_2}A_3\xrightarrow{f_3}\cdots\xrightarrow{f_{d}}A_{d+1}\xrightarrow{f_{d+1}}\Sigma^d A_0\xrightarrow{(-1)^{d}\Sigma^d f_0}\Sigma^d A_1.$$
A \emph{morphism} of $(d+2)$-$\Sigma^d$-sequences is  a sequence of morphisms $\varphi=(\varphi_0,\varphi_1,\cdots,\varphi_{d+1})$ such that the following diagram
$$\xymatrix{
A_0 \ar[r]^{f_0}\ar[d]^{\varphi_0} & A_1 \ar[r]^{f_1}\ar[d]^{\varphi_1} & A_2 \ar[r]^{f_2}\ar[d]^{\varphi_2} & \cdots \ar[r]^{f_{d}}& A_{d+1} \ar[r]^{f_{d+1}}\ar[d]^{\varphi_{d+1}} & \Sigma^d A_0 \ar[d]^{\Sigma^d \varphi_0}\\
B_0 \ar[r]^{g_0} & B_1 \ar[r]^{g_1} & B_2 \ar[r]^{g_2} & \cdots \ar[r]^{g_{d}}& B_{d+1} \ar[r]^{g_{d+1}}& \Sigma^d B_0
}$$
commutes, where each row is a $(d+2)$-$\Sigma^d$-sequence. It is an {\em isomorphism} if $\varphi_0, \varphi_1, \varphi_2, \cdots, \varphi_{d+1}$ are all isomorphisms in $\C$, and
a \emph{weak isomorphism} if $\varphi_i$ and $\varphi_{i+1}$ are isomorphisms for some $0\leq i\leq d+1$ (with
$\varphi_{d+2}:=\Sigma\varphi_0$). Note that the composition of two weak isomorphisms need not be a
weak isomorphism.
\medskip

We recall the notion of a right $(d+2)$-angulated category from \cite[Definition 2.1]{l2}.
The axiom (N4) presented below is based on the higher octahedral axiom
introduced in \cite{ahbt}.

\begin{definition}\label{d1}
A {\em right} $(d+2)$-\emph{angulated category} is a triple $(\C, \Sigma^d, \Theta)$, where $\C$ is an additive category, $\Sigma^d$ is the $d$-suspension functor of $\C$, and $\Theta$ is a class of $(d+2)$-$\Sigma^d$-sequences (whose elements are called right $(d+2)$-angles), which satisfies the following axioms:
\begin{itemize}
\item[(N1)]
\begin{itemize}
\item[(a)] The class $\Theta$ is closed under isomorphisms, direct sums and  direct summands.

\item[(b)] For each object $A\in\C$ the trivial sequence
$$A\xrightarrow{{\rm Id}_A} A\rightarrow 0\rightarrow 0\rightarrow\cdots\rightarrow 0\rightarrow \Sigma^dA$$
belongs to $\Theta$.

\item[(c)] Each morphism $f_0\colon A_0\rightarrow A_1$ in $\C$ can be extended to a right $(d+2)$-angle: $$A_0\xrightarrow{f_0}A_1\xrightarrow{f_1}A_2\xrightarrow{f_2}\cdots\xrightarrow{f_{d-1}}A_d\xrightarrow{f_d}A_{d+1}\xrightarrow{f_{d+1}}\Sigma^d A_0.$$
\end{itemize}
\item[(N2)] If a $(d+2)$-$\Sigma^d$-sequence belongs to $\Theta$, then its left rotation belongs to $\Theta$.

\item[(N3)] Each solid commutative diagram
$$\xymatrix{
A_0 \ar[r]^{f_0}\ar[d]^{\varphi_0} & A_1 \ar[r]^{f_1}\ar[d]^{\varphi_1} & A_2 \ar[r]^{f_2}\ar@{-->}[d]^{\varphi_2} & \cdots \ar[r]^{f_{d}}& A_{d+1} \ar[r]^{f_{d+1}}\ar@{-->}[d]^{\varphi_{d+1}} & \Sigma^d A_0 \ar[d]^{\Sigma^d \varphi_0}\\
B_0 \ar[r]^{g_0} & B_1 \ar[r]^{g_1} & B_2 \ar[r]^{g_2} & \cdots \ar[r]^{g_{d}}& B_{d+1} \ar[r]^{g_{d+1}}& \Sigma^d B_0
}$$ with rows in $\Theta$ can be completed to a morphism of  $(d+2)$-$\Sigma^d$-sequences.

\item[(N4)]  Given the solid part of the diagram
  \begin{center}
    \begin{tikzpicture}
      \diagram{d}{2.5em}{2.5em}{
        A_0 & A_1 & A_2 & A_3 & \cdots & A_{n} & A_{n+1} & \Sigma A_0\\
        A_0 & B_1 & B_2 & B_3 & \cdots & B_{n} & B_{n+1} & \Sigma
        A_0\\
        A_1 & B_1 & C_2 & C_3 & \cdots & C_{n} & C_{n+1} & \Sigma
        A_1\\
      };

      \path[->,midway,font=\scriptsize]
        (d-1-1) edge node[above] {$f_0$} (d-1-2)
        ([xshift=-0.1em] d-1-1.south) edge[-] ([xshift=-0.1em] d-2-1.north)
        ([xshift=0.1em] d-1-1.south) edge[-] ([xshift=0.1em] d-2-1.north)
        (d-1-2) edge node[above] {$f_1$} (d-1-3)
                     edge node[right] {$\varphi_1$} (d-2-2)
        (d-1-3) edge node[above] {$f_2$} (d-1-4)
                     edge[densely dashed] node[right] {$\varphi_2$} (d-2-3)
        (d-1-4) edge node[above] {$f_3$} (d-1-5)
                     edge[densely dashed] node[right] {$\varphi_3$} (d-2-4)
        (d-1-5) edge node[above] {$f_{n-1}$} (d-1-6)
        (d-1-6) edge node[above] {$f_{n}$} (d-1-7)
                     edge[densely dashed] node[right] {$\varphi_{n}$} (d-2-6)
        (d-1-7) edge node[above] {$f_{n+1}$} (d-1-8)
                     edge[densely dashed] node[right] {$\varphi_{n+1}$} (d-2-7)
        ([xshift=-0.1em] d-1-8.south) edge[-] ([xshift=-0.1em] d-2-8.north)
        ([xshift=0.1em] d-1-8.south) edge[-] ([xshift=0.1em] d-2-8.north)
        (d-2-1) edge node[above] {$g_0$} (d-2-2)
                     edge node[right] {$f_0$} (d-3-1)
        (d-2-2) edge node[above] {$g_1$} (d-2-3)
        ([xshift=-0.1em] d-2-2.south) edge[-] ([xshift=-0.1em] d-3-2.north)
        ([xshift=0.1em] d-2-2.south) edge[-] ([xshift=0.1em] d-3-2.north)
        (d-2-3) edge node[above] {$g_2$} (d-2-4)
                     edge[densely dashed] node[right] {$\theta_2$} (d-3-3)
        (d-2-4) edge node[above] {$g_3$} (d-2-5)
                     edge[densely dashed] node[right] {$\theta_3$} (d-3-4)
        (d-2-5) edge node[above] {$g_{n - 1}$} (d-2-6)
        (d-2-6) edge node[above] {$g_{n}$} (d-2-7)
                     edge[densely dashed] node[right] {$\theta_{n}$} (d-3-6)
        (d-2-7) edge node[above] {$g_{n+1}$} (d-2-8)
                     edge[densely dashed] node[right] {$\theta_{n+1}$} (d-3-7)
        (d-2-8) edge node[right] {$\Sigma f_0$} (d-3-8)
        (d-3-1) edge node[above] {$\varphi_1$} (d-3-2)
        (d-3-2) edge node[above] {$h_1$} (d-3-3)
        (d-3-3) edge node[above] {$h_2$} (d-3-4)
        (d-3-4) edge node[above] {$h_3$} (d-3-5)
        (d-3-5) edge node[above] {$h_{n - 1}$} (d-3-6)
        (d-3-6) edge node[above] {$h_{n}$} (d-3-7)
        (d-3-7) edge node[above] {$h_{n+1}$} (d-3-8)
        (d-1-4) edge[densely dashed,out=-102,in=30] node[pos=0.15,left] {$\psi_3$}
          (d-3-3)
        (d-1-7) edge[densely dashed,out=-102,in=30] node[pos=0.15,left] {$\psi_{n+1}$}
          (d-3-6);
    \end{tikzpicture}
  \end{center}
  with commuting squares and with rows in $\Theta$, the dotted
  morphisms exist such that each square commutes,  and the
  $(d+2)$-$\Sigma^d$-sequence
$$A_2\xrightarrow{\left(
                    \begin{smallmatrix}
                      f_2 \\
                      \varphi_2 \\
                    \end{smallmatrix}
                  \right)} A_3\oplus B_2\xrightarrow{\left(
                             \begin{smallmatrix}
                               -f_3 & 0 \\
                               \varphi_3 & -g_2 \\
                               \phi_3 & \psi_2 \\
                             \end{smallmatrix}
                           \right)}
 A_4\oplus B_3\oplus C_2\xrightarrow{\left(
                                       \begin{smallmatrix}
                                         -f_4 & 0 & 0 \\
                                         -\varphi_4 & -g_3 & 0 \\
                                         \phi_4 & \psi_3 & h_2 \\
                                       \end{smallmatrix}
                                     \right)}A_5\oplus B_4\oplus C_3$$
$$\xrightarrow{\left(
                                       \begin{smallmatrix}
                                         -f_5 & 0 & 0 \\
                                         \varphi_5 & -g_4 & 0 \\
                                         \phi_5 & \psi_4 & h_3 \\
                                       \end{smallmatrix}
                                     \right)}\cdots\xrightarrow{\scriptsize\left(\begin{smallmatrix}
             -f_{d} & 0 & 0 \\
             (-1)^{d+1}\varphi_{d} & -g_{d-1} & 0 \\
             \phi_{d} & \psi_{d-1} & h_{d-2} \\
             \end{smallmatrix}
             \right)}A_{d+1}\oplus B_{d}\oplus C_{d-1}$$
$$\xrightarrow{\left(
                                                       \begin{smallmatrix}
                                                         (-1)^{d+1}\varphi_{d+1} &-g_{d} &0 \\
                                                          \phi_{d+1}& \psi_{d}& h_{d-1} \\
                                                       \end{smallmatrix}
                                                     \right)}B_{d+1}\oplus C_{d}\xrightarrow{(\psi_{d+1},\ h_{d})}C_{d+1}\xrightarrow{\Sigma^d f_1\circ h_{d+1}}\Sigma^d A_2 \hspace{10mm}$$
belongs to $\Theta$.
   \end{itemize}
\end{definition}
The notion of a \emph{left $(d+2)$-angulated category} is defined dually.
\vspace{1mm}

If $\Sigma^d$ is an automorphism, it is easy to see that the converse of an axiom (N2) also holds, thus the right
$(d+2)$-angulated category $(\C,\Sigma^d, \Theta)$ is a $(d+2)$-angulated category in the sense of Geiss-Keller-Oppermann \cite[Definition 1.1]{gko} and in the sense of Bergh-Thaule \cite[Theorem 4.4]{bt1}. If $(\C,\Sigma^d, \Theta)$ is a right $(d+2)$-angulated category, $(\C,\Omega, \Phi)$
is a left $(d+2)$-angulated category, $\Omega$ is a quasi-inverse of $\Sigma^d$ and $\Theta=\Phi$, then $(\C,\Sigma^d, \Theta)$ is a $(d+2)$-angulated category.
\medskip

Suppose that $\Sigma^d$ is an automorphism. Now modify axiom (N1) into a new axiom (N1$^\ast$).

\begin{itemize}
\item[(N1)$^\ast$]
\begin{itemize}
\item[(a)] The class $\Theta$ is closed under weak isomorphisms.

\item[(b)] For each object $A\in\C$ the trivial sequence
$$A\xrightarrow{{\rm Id}_A} A\rightarrow 0\rightarrow 0\rightarrow\cdots\rightarrow 0\rightarrow \Sigma^dA$$
belongs to $\Theta$.

\item[(c)] Each morphism $f_0\colon A_0\rightarrow A_1$ in $\C$ can be extended to a $(d+2)$-angle: $$A_0\xrightarrow{f_0}A_1\xrightarrow{f_1}A_2\xrightarrow{f_2}\cdots\xrightarrow{f_{d-1}}A_d\xrightarrow{f_d}A_{d+1}\xrightarrow{f_{d+1}}\Sigma^d A_0.$$
\end{itemize}
\end{itemize}

\begin{remark}\cite[Theorem 3.2]{bt1}
If $\Theta$ is a collection of $(d+2)$-$\Sigma^d$-sequences satisfying the axioms (N2) and
(N3), and $\Sigma^d$ is an automorphism. Then the following are equivalent:
\begin{itemize}
\item[\rm (1)] $\Theta$  satisfies (N1);
\item[\rm (2)] $\Theta$ satisfies (N1)$^{\ast}$.
\end{itemize}
\end{remark}

Let $\C$ be an additive category and $f\colon A\rightarrow B$ a morphism in $\C$. A \emph{weak cokernel} of $f$ is a morphism
$g\colon B\rightarrow C$ such that for any $X\in\C$ the sequence of abelian groups
$$\Hom_\C(C,X)\xrightarrow{~g^\ast~}\Hom_\C(B,X)\xrightarrow{~f^\ast~}\Hom_\C(A,X)$$
is exact. Equivalently, $g$ is a weak cokernel of $f$ if $gf=0$ and for each morphism
$h\colon B\rightarrow X$ such that $hf=0$ there exists a (not necessarily unique) morphism
$p\colon C\rightarrow X$ such that $h=pg$. These properties are subsumed in the following
commutative diagram
$$\xymatrix{A\ar[r]^{f}\ar[dr]_{0} & B\ar[r]^{g}\ar[d]^{h}&C\ar@{-->}[dl]^{p}  \\
& X &}$$
Clearly, a weak cokernel $g$ of $f$ is a cokernel of $f$ if and only if $g$ is an epimorphism.
The concept of a \emph{weak kernel} is defined dually.

\begin{lemma}\label{bu1}
Let $(\C, \Sigma^d, \Theta)$ be a $(d+2)$-angulated category and
$$A_0\xrightarrow{f_0}A_1\xrightarrow{f_1}A_2\xrightarrow{f_2}\cdots\xrightarrow{f_{d-1}}A_d\xrightarrow{f_d}A_{d+1}\xrightarrow{f_{d+1}}\Sigma^d A_0.$$
be a $(d+2)$-angle. Then
$f_{i}\circ f_{i-1}=0, ~i=1,2,\cdots,d+1$. That is to say, any composition of
two consecutive morphisms in a $(d+2)$-angle vanishes.
\end{lemma}

\proof  This is \cite[Prposition 2.5 (a)]{gko}. See also \cite[Lemma 3.1]{bt1}.  \qed

\begin{lemma}\label{bu2}
Let $(\C, \Sigma^d, \Theta)$ be a $(d+2)$-angulated category and
$$A_0\xrightarrow{f_0}A_1\xrightarrow{f_1}A_2\xrightarrow{f_2}\cdots\xrightarrow{f_{d-1}}A_d\xrightarrow{f_d}A_{d+1}\xrightarrow{f_{d+1}}\Sigma^d A_0.$$
be a $(d+2)$-angle. Then $f_{i}$ is a weak cokernel of $f_{i-1}$
and $f_{i-1}$ is a weak kernel of $f_{i}$, $i=1,2,\cdots,d+1$.
\end{lemma}

\proof  This follows from \cite[Prposition 2.5 (a) and Remark 2.2 (c)]{gko}.
See also \cite[Lemma 2.4]{l1}. \qed

\subsection{Auslander-Reiten $(d+2)$-angles}
\setcounter{equation}{0}
Let $\C$ be a $(d+2)$-angulated category. We denote by ${\rm rad}_{\C}$ the Jacobson radical of $\C$. Namely, it is given by the formula
$${\rm rad}_{\C}(X,Y)=\{g\colon X\to Y~|~{\rm Id_{\emph{X}}-}hg~\mbox{is invertible for any}~ h\colon Y\to X\},$$
for all objects $X$ and $Y$ in $\C$.

Next, we recall from \cite{ar} some terminology for the Auslander--Reiten
theory. Let $f\colon X \to Y$ be a morphism. One says that $f$ is \emph{right almost split} if $f$ is not a split epimorphism and every non-split epimorphism  $g\colon M \to  Y$ factors through $f$.
In a dual manner, one defines $f$ to be \emph{left almost split}.

\begin{definition}\cite[Definition 3.8]{iy}, \cite[Definition 5.1]{f}
Let $\C$ be a $(d+2)$-angulated category.
A $(d+2)$-angle
$$A_{\bullet}:~~\xymatrix {A_0 \xrightarrow{~\alpha_0~}A_1 \xrightarrow{~\alpha_1~} A_2 \xrightarrow{~\alpha_2~} \cdots
  \xrightarrow{~\alpha_{d - 1}~} A_d \xrightarrow{~\alpha_{d}~} A_{d+1}\xrightarrow{~\alpha_{d+1}~} \Sigma^d A_0}$$
in $\C$ is called an \emph{Auslander-Reiten $(d+2)$-angle }if
$\alpha_0$ is left almost split, $\alpha_d$ is right almost split and
when $d\geq 2$, also $\alpha_1,\alpha_2,\cdots,\alpha_{d-1}$ are in $\rad_{\C}$.
\end{definition}

\begin{remark}\cite[Remark 5.2]{f}\label{rem}
Assume $A_{\bullet}$ as in the above definition is an Auslander-Reiten $(d+2)$-angle.
Since $\alpha_0$ is left almost split which implies that $\End(A_0)$ is local and hence $A_0$ is indecomposable. Similarly, since $\alpha_d$ is right almost split, then $\End(A_{d+1})$ is local and hence $A_{d+1}$ is indecomposable.
Moreover, when $d=1$, we have $\alpha_0$ and $\alpha_d$ in $\rad_{\C}$, so that $\alpha_d$ is right minimal and $\alpha_0$ is left minimal. When $d\geq 2$, since $\alpha_{d-1}\in\rad_{\C}$, we have that $\alpha_d$ is right minimal and similarly $\alpha_0$ is left
minimal.
\end{remark}

\begin{remark}\cite[Lemma 5.3]{f}
Let $\C$ be a $(d+2)$-angulated category and
$$A_{\bullet}:~~\xymatrix {A_0 \xrightarrow{~\alpha_0~}A_1 \xrightarrow{~\alpha_1~} A_2 \xrightarrow{~\alpha_2~} \cdots
  \xrightarrow{~\alpha_{d - 1}~} A_d \xrightarrow{~\alpha_{d}~} A_{d+1}\xrightarrow{~\alpha_{d+1}~} \Sigma^d A_0}$$
be a $(d+2)$-angle in $\C$. Then the following statements are equivalent:
\begin{itemize}
\item[(1)] $A_{\bullet}$ is an Auslander-Reiten $(d+2)$-angle;
\item[(2)] $\alpha_0,\alpha_1,\cdots,\alpha_{d-1}$ are in ${\rm rad}_{\C}$ and $\alpha_d$ is right almost split;
\item[(3)] $\alpha_1,\alpha_2,\cdots,\alpha_{d}$ are in ${\rm rad}_{\C}$ and $\alpha_0$ is left almost split.
\end{itemize}
\end{remark}

\begin{lemma}\emph{\cite[Lemma 5.4]{f}}\label{keylemma}
Let $\C$ be a $(d+2)$-angulated category and
$$A_{\bullet}:~~\xymatrix {A_0 \xrightarrow{~\alpha_0~}A_1 \xrightarrow{~\alpha_1~} A_2 \xrightarrow{~\alpha_2~} \cdots
  \xrightarrow{~\alpha_{d - 1}~} A_d \xrightarrow{~\alpha_{d}~} A_{d+1}\xrightarrow{~\alpha_{d+1}~} \Sigma^d A_0}$$
be a $(d+2)$-angle in $\C$.
Assume that $\alpha_d$ is right almost split and if $d\geq2$, also that $\alpha_1,\alpha_2,\cdots,\alpha_{d-1}$ are in $\emph{\rad}_{\C}$.
Then the following are equivalent:
\begin{itemize}
\item[\emph{(1)}] $A_{\bullet}$ is an Auslander-Reiten $(d+2)$-angle;
\item[\emph{(2)}] $\emph{\End}(A_0)$ is local;
\item[\emph{(3)}] $\alpha_{d+1}$ is left minimal;
\item[\emph{(4)}] $\alpha_{0}$ is in $\emph{\rad}_{\C}$.
\end{itemize}
\end{lemma}
\medskip

In the case $d=1$, so in the case of a triangulated category, a morphism can be extended to a triangle
in a unique way up to isomorphism. On the other hand, for $d\geq 2$, a morphism can be extendend to a
$(d+2)$-angle in different non-isomorphic ways. However, we still have a unique ``minimal" $(d+2)$-angle
extending any given morphism.

\begin{lemma}\emph{\cite[Lemma 5.18]{ot}, \cite[Lemma 3.14]{f}}\label{lemma}
Let $d\geq 2$ and $h\colon A_{d+1}\to \Sigma^d A_0$ be any morphism in a $(d+2)$-angulated category $\C$. Then, up to
isomorphism, there exists a unique $(d+2)$-angle of the form
$$\xymatrix {A_0 \xrightarrow{~\alpha_0~}A_1 \xrightarrow{~\alpha_1~} A_2 \xrightarrow{~\alpha_2~} \cdots
  \xrightarrow{~\alpha_{d-2}~} A_{d-1} \xrightarrow{~\alpha_{d - 1}~} A_d \xrightarrow{~\alpha_{d}~} A_{d+1}\xrightarrow{~h~} \Sigma^d A_0}$$
with $\alpha_1,\alpha_2,\cdots,\alpha_{d-1}$ in ${\rm rad_{\C}}$.
\end{lemma}

\section{Every Serre functor is a $(d+2)$-angulated functor}
\setcounter{equation}{0}

Let $\C$ be a $k$-linear
Hom-finite additive category.
When there is no danger of confusion, we will sometimes instead of
$$\Hom_\C(X,A)\xrightarrow{~\Hom_\C(X,\ f)~}\Hom_\C(X,B)$$
 write one of the following simplified forms:
$$\Hom(X,A)\xrightarrow{~\Hom(X,\ f)~}\Hom(X,B)$$
$$(X,A)\xrightarrow{~(X,\ f)~}(X,B).$$

A \emph{{{right Serre}} functor}
is a $k$-linear additive functor $\F\colon\C\to \C$ together with isomorphisms
$$\eta_{A,B}\colon\Hom(A,B)\to \Hom(B,\F A)^\ast$$
for any $A,B\in\C$ which are natural in $A$ and $B$, where $(-)^{\ast}:=\Hom_{k}(-,k)$.

 A \emph{left Serre functor} is a $k$-linear additive functor
$\F\colon\C\to\C$ together with  isomorphisms
$$\zeta_{A,B}\colon \Hom(A,B)\to \Hom(\F B,A)^\ast$$
for any $A,B\in\C$ which are natural in $A$ and $B$.

Let $\eta_A\colon\Hom(A,\F A)\to k$ be given by $\eta_{A,A}({\rm Id}_A)$ and
let $f\in\Hom(A, B)$. Looking at the following commutative diagram
(which follows from the naturality of $\eta_{A,B}$ in $B$)
$$\xymatrix@C=1.2cm @R=1.2cm{\Hom(A,A)\ar[r]^{\eta_{A,A}\quad}\ar[d]_{\Hom(A,\ f)}&\Hom(A,\F A)^\ast\ar[d]^{\Hom(f, \ \F A)^\ast}\\
\Hom(A,B)\ar[r]^{\eta_{A,B}\quad}& \Hom(B,\F A)^\ast}$$
we find for $g\in\Hom(B,\F A)$\
$$\eta_{A,B}(f)(g)=\eta_A (g f).$$
Similarly by the naturality of $\eta_{A,B}$ in $A$ we obtain a commutative diagram
$$\xymatrix@C=1.2cm @R=1.2cm{\Hom(B,B)\ar[r]^{\eta_{B,B}\quad}\ar[d]_{\Hom(f,\ B)}&\Hom(B,\F B)^\ast\ar[d]^{\Hom(B,\ \F f)^\ast}\\
\Hom(A,B)\ar[r]^{\eta_{A,B}\quad}& \Hom(B,\F A)^\ast}$$
This yields for $g\in \Hom(B,\F A)$ the formula
\begin{equation}\label{t1}
\begin{array}{l}
\eta_{A,B}(f)(g)=\eta_B(\F(f) g).
\end{array}
\end{equation}

\begin{remark}\cite[Lemma I.1.5]{rv}
 $\C$ has a Serre functor if and only it has both a
  right and a left Serre functor if and only it has a right Serre functor which is an auto-equivalence.
\end{remark}

\begin{definition}\cite{ja,bt2}
Let $\C$ be a $(d+2)$-angulated category with $d$-suspension functor $\Sigma^d$. An (covariant) additive functor $\F\colon \C\to \C$ is called \emph{$(d+2)$-angulated} if it has the following properties:
\begin{itemize}

\item[$(1)$] There exists a natural isomorphism $\phi\colon \F\Sigma^d \to \Sigma^d\F$;
\item[$(2)$] $\F$ preserves $(d+2)$-angles, that is to say, if
$$A_0\xrightarrow{f_0}A_1\xrightarrow{f_1}A_2\xrightarrow{f_2}\cdots\xrightarrow{f_d}A_{d+1}\xrightarrow{f_{d+1}}\Sigma^d A_0$$
is a $(d+2)$-angle in $\C$, then
$$\F A_0\xrightarrow{\F(f_0)}\F A_1\xrightarrow{\F(f_1)}\F A_2\xrightarrow{\F(f_2)}\cdots\xrightarrow{\F(f_d)}\F A_{d+1}\xrightarrow{\phi_{A_0}\circ \F(d_{d+1})}\Sigma^d\F A_0$$
is a $(d+2)$-angle in $\C$.
\end{itemize}
\end{definition}

Bondal and Kapranov \cite[Proposition 3.3]{bk} prove that every Serre functor is a triangulated functor in a triangulated category.
Van den Bergh has since provided a different proof, see \cite[Theorem A.4.4]{b}.
The following result  shows that the Serre functor is a $(d+2)$-angulated functor in a $(d+2)$-angulated category. This generalizes the work by Bondal and Kapranov.
Our proof is an adaptation of the proof of Van den Bergh. For more details, see also \cite[Theorem 10.5.1]{z}.

\begin{theorem}\label{main1}
Let $\C$ be a $(d+2)$-angulated category with right Serre functor $\F$.
Then there exists a  natural isomorphism $\zeta\colon \F\Sigma^d\to \Sigma^d\F$ such that
$\F\colon \C\to\C$ is a $(d+2)$-angulated functor.
\end{theorem}

\proof Assume that $(\F,\eta_{X,Y})$ is a right Serre functor. Put
$$\eta_{X}:=\eta_{X,X}({\rm Id}_X)\in\Hom(X,\F X)^{\ast},~\mbox{for all}~X\in\C.$$

\textbf{Step 1:} We claim that there exists a  natural isomorphism $\zeta\colon \F\Sigma^d\to \Sigma^d\F$ such that
\begin{equation}\label{t2}
\begin{array}{l}
\eta_{X}(\Sigma^{-d}\zeta_X\circ\Sigma^{-d}f)=(-1)^d\eta_{\Sigma^dX}(f),~\mbox{for all}~ f\in\Hom(\Sigma^dX,\F(\Sigma^dX)) ~\textrm{and}~ X\in\C.
\end{array}
\end{equation}
By the following two isomorphisms
$$\Hom(\Sigma^dX,\F(\Sigma^dX))^{\ast}\xrightarrow{~~\sim~~}\Hom(X,\Sigma^{-d}\F(\Sigma^dX))^{\ast}$$
$$\eta^{\ast}_{X,~\Sigma^{-d}\F(\Sigma^dX)}\colon\Hom(X,\Sigma^{-d}\F(\Sigma^dX))^{\ast}
\xleftarrow{~~\sim~~}\Hom(\Sigma^{-d}\F(\Sigma^dX),\F X)$$
and $(-1)^d\eta_{\Sigma^dX}\in\Hom(\Sigma^dX,\F(\Sigma^dX))^{\ast}$, there exists
$\zeta_X\in\Hom(\F(\Sigma^dX),\Sigma^d(\F X))$ such that
\begin{equation}\label{t3}
\begin{array}{l}
\eta^{\ast}_{X,~\F\Sigma^{-d}(\Sigma^dX)}(\Sigma^{-d}\zeta_X)(g)=(-1)^d\eta_{\Sigma^dX}(\Sigma^dg),
~\mbox{for all}~ g\colon
X\to\F\Sigma^{-d}(\Sigma^dX).
\end{array}
\end{equation}
Using Serre duality, the above equality (\ref{t3}) is equivalent to
\begin{equation}\label{t4}
\begin{array}{l}
\eta_{X,~\Sigma^{-d}\F(\Sigma^dX)}(g)(\Sigma^{-d}\zeta_X)=(-1)^d\eta_{\Sigma^dX}(\Sigma^dg).
\end{array}
\end{equation}
By the equality (\ref{t1}),  the above equality (\ref{t4}) is identified with
$$\eta_{X}(\Sigma^{-d}\zeta_X\circ \Sigma^{-d}f)=(-1)^d\eta_{\Sigma^dX}(f),
~\mbox{for all}~f\in\Hom(\Sigma^dX,\F(\Sigma^dX)),$$
as desired. That is to say, the equality (\ref{t2}) holds.

It remains to show that $\zeta$ is a natural isomorphism.

By the following two isomorphisms
$$\eta^{\ast}_{\Sigma^dX,~\Sigma^d(\F X)}\colon\Hom(\Sigma^d(\F X),\F(\Sigma^{d}X))
\xrightarrow{~~\sim~~}\Hom(\Sigma^{d}X,\Sigma^d(\F X))^{\ast}$$
$$\Hom(\Sigma^dX,\Sigma^d(\F X))^{\ast}\xleftarrow{~~\sim~~}\Hom(X,\F X)^{\ast}$$
and $(-1)^d\eta_{X}\in\Hom(X,\F X)^{\ast}$, there exists
$\theta_X\in\Hom(\Sigma^d(\F X),\F(\Sigma^dX))$ such that
\begin{equation}\label{t5}
\begin{array}{l}
\eta^{\ast}_{\Sigma^dX,~\Sigma^d(\F X)}(\theta_X)(g)=(-1)^d\eta_{X}(\Sigma^{-d}g),~\mbox{for all}~g\colon
\Sigma^dX\to\Sigma^d(\F X).
\end{array}
\end{equation}
Using Serre duality, the above equality (\ref{t5}) is equivalent to
\begin{equation}\label{t6}
\begin{array}{l}
\eta^{\ast}_{\Sigma^dX,~\Sigma^d(\F X)}(g)(\theta_X)=(-1)^d\eta_X(\Sigma^{-d}g).
\end{array}
\end{equation}
By the equality (\ref{t1}),  the above equality (\ref{t6}) is identified with
\begin{equation}\label{t7}
\begin{array}{l}
\eta_{\Sigma^dX}(\theta_X\Sigma^dg)=(-1)^d\eta_{X}(g),
~\mbox{for all}~g\in\Hom(X,\F X).
\end{array}
\end{equation}
We now show that $\zeta$ is an isomorphism.  For any morphism $g\colon X\to \F X$, we have
$$\eta_X(\Sigma^{-d}(\zeta_X\theta_X)g)\xlongequal{(\ref{t2})}(-1)^d\eta_{\Sigma^dX}(\theta_X\Sigma^dg)
\xlongequal{(\ref{t7})}\eta_X(g).$$
By the bilinear map
$(-,-)\colon\Hom(X,\F X)\times \Hom(\F X,\F X)\longrightarrow k$
is non-degenerate, we obtain $\Sigma^{-d}(\zeta_X\theta_X)={\rm Id}_{\F X}$ and then
$\zeta_X\theta_X={\rm Id}_{\Sigma^d(\F X)}$.

Similarly, for any morphism $g\colon \Sigma^dX\to \F(\Sigma^dX)$, we have
$$\eta_{\Sigma^dX}(\theta_X\zeta_Xg)\xlongequal{(\ref{t7})}(-1)^d\eta_{X}(\Sigma^{-d}(\zeta_Xg))
\xlongequal{(\ref{t2})}\eta_{\Sigma^dX}(g).$$
So $\theta_X\zeta_X={\rm Id}_{\F(\Sigma^dX)}$.
This shows that $\zeta_X$ is an isomorphism.

We now show that $\zeta$ is a natural transformation. Namely, for any morphism $f\colon X\to Y$, we have the following commutative diagram:
$$\xymatrix@C=1.2cm @R=1.2cm{\F(\Sigma^dX)\ar[r]^{\zeta_X}\ar[d]_{\F(\Sigma^df)}&\Sigma^d(\F X)\ar[d]^{\Sigma^d(\F f)}\\
\F(\Sigma^dY)\ar[r]^{\zeta_Y}& \Sigma^d(\F Y).}$$
Note that $\Sigma^{-d}\zeta_X\in\Hom(\Sigma^{-d}\F(\Sigma^dX),\F X)$ and
$\Sigma^{-d}\zeta_Y\in\Hom(\Sigma^{-d}\F(\Sigma^dY),\F Y )$. For any morphism
$h\in\Hom(Y,\Sigma^{-d}\F(\Sigma^dX))$, we have
$$\begin{aligned}%
\Hom(f,-)^{\ast}(\eta^{\ast}_{X,~\Sigma^{-d}\F(\Sigma^dX)}(\Sigma^{-d}\zeta_X))(h)&\xlongequal{\rm duality}
\eta^{\ast}_{X,~\Sigma^{-d}\F(\Sigma^dX)}(\Sigma^{-d}\zeta_X)(hf)\\
                     &\xlongequal{\rm duality}\eta_{X,~\Sigma^{-d}\F(\Sigma^dX)}(hf)(\Sigma^{-d}\zeta_X)\\
                     &\xlongequal{(\ref{t1})}\eta_X(\Sigma^{-d}\zeta hf)\\
                     &\xlongequal{(\ref{t2})}(-1)^d\eta_{\Sigma^dX}(\Sigma^dh\circ\Sigma^df)\\
                     &\xlongequal{(\ref{t1})}(-1)^d\eta_{\Sigma^dY}(\F(\Sigma^df)\Sigma^dh)\\
                     &\xlongequal{(\ref{t2})}(-1)^d\eta_{Y}(\Sigma^{-d}\zeta_Y\Sigma^{-d}\F(\Sigma^df)h)\\
                     &\xlongequal{(\ref{t1})}(-1)^d\eta_{Y,~\Sigma^{-d}\F(\Sigma^dY)}(\Sigma^{-d}\F(\Sigma^df)h)(\Sigma^{-d}\zeta_Y)\\
                     &\xlongequal{\rm duality}(-1)^d\eta^{\ast}_{Y,~\Sigma^{-d}\F(\Sigma^dY)}(\Sigma^{-d}\zeta_Y)(\Sigma^{-d}\F(\Sigma^df)h)\\
                     &\xlongequal{\rm duality}\Hom(-,\Sigma^{-d}\F(\Sigma^df))^{\ast}(\eta^{\ast}_{Y,~\Sigma^{-d}\F(\Sigma^dY)}(\Sigma^{-d}\zeta_Y))(h).
\end{aligned}$$
Consider the following commutative diagram:
$$\xymatrix@C=2cm @R=1cm{(\Sigma^{-d}\F(\Sigma^dX),~\F X)\ar[r]^{\eta^{\ast}_{X,~\Sigma^{-d}\F(\Sigma^dX)}}\ar[d]_{(-,~\F f)}&(X,~\Sigma^{-d}\F(\Sigma^dX))^{\ast}\ar[d]^{(f,~-)^{\ast}}\\
(\Sigma^{-d}\F(\Sigma^dX),~\F Y)\ar[r]^{\eta^{\ast}_{Y,~\Sigma^{-d}\F(\Sigma^dX)}}& (Y,~\Sigma^{-d}\F(\Sigma^dX))^{\ast}\\
(\Sigma^{-d}\F(\Sigma^dY),~\F Y)\ar[r]^{\eta^{\ast}_{Y,~\F\Sigma^{-d}(\Sigma^dY)}}\ar[u]^{(\Sigma^{-d}\F(\Sigma^df),~-)}&(Y,~\Sigma^{-d}\F(\Sigma^dY))^{\ast}
\ar[u]_{(-,~\Sigma^{-d}\F(\Sigma^df))^{\ast}}}$$
It follows that
$\Hom(-,\F f)(\Sigma^{-d}\zeta_X)=\Hom(\Sigma^{-d}\F(\Sigma^df),~-))(\Sigma^{-d}\zeta_Y)$
and then
$$\F(f)\Sigma^{-d}\zeta_X=\Sigma^{-d}\zeta_Y\Sigma^{-d}\F(\Sigma^df).$$
Hence $\Sigma^{d}(\F f)\zeta_X=\zeta_Y\F(\Sigma^df)$ as required.

\textbf{Step 2:}  Assume that
$$A_0\xrightarrow{\alpha_0}A_1\xrightarrow{\alpha_1}A_2\xrightarrow{\alpha_2}\cdots
\xrightarrow{\alpha_{d-1}}A_d\xrightarrow{\alpha_d}A_{d+1}\xrightarrow{\alpha_{d+1}}\Sigma^dA_{0}$$
is any $(d+2)$-angle in $\C$.  The morphism $\F\alpha_0$ can be embedded in a $(d+2)$-angle
$$\F A_0\xrightarrow{\F\alpha_0}\F A_1\xrightarrow{\beta_1}B_2\xrightarrow{\beta_2}\cdots
\xrightarrow{\beta_{d-1}}B_d\xrightarrow{\beta_d}B_{d+1}\xrightarrow{\beta_{d+1}}\Sigma^d\F A_{0}.$$
In order to prove that $(\F,\zeta)$ is a $(d+2)$-angulated functor, we only need to show that there
are morphisms $\omega_i\colon B_i\to\F A_i,~i=2,3,\cdots,d+1$ such that the following diagram commutes:
$$\xymatrix@C=1.2cm@R=1.2cm{\F A_0\ar[r]^{\F\alpha_0}\ar@{=}[d]&\F A_1\ar[r]^{\beta_1}\ar@{=}[d]&B_2\ar[r]^{\beta_2}\ar@{-->}[d]^{\omega_2}&B_3\ar[r]^{\beta_3}\ar@{-->}[d]^{\omega_3}
&\cdot\cdot\cdot\ar[r]^{\beta_{d-1}}
&B_d\ar[r]^{\beta_d}\ar@{-->}[d]^{\omega_{d}}&B_{d+1}\ar[r]^{\beta_{d+1}}\ar@{-->}[d]^{\omega_{d+1}}&\Sigma^d\F A_0\ar@{=}[d]\\
\F A_0\ar[r]^{\F\alpha_0}&\F A_1\ar[r]^{\F\alpha_1}&\F A_2\ar[r]^{\F\alpha_2}&\F A_3\ar[r]^{\F\alpha_3}&\cdot\cdot\cdot
\ar[r]^{\F\alpha_{d-1}}&\F A_d\ar[r]^{\F\alpha_d}&\F A_{d+1}\ar[r]^{\zeta_{A_0}\F\alpha_{d+1}}&\Sigma^d\F A_0.}$$
Note that the class of $(d+2)$-angles is closed under weak isomorphisms, then
$$\F A_0\xrightarrow{\F\alpha_0}\F A_1\xrightarrow{\F\alpha_1}\F A_2\xrightarrow{\F\alpha_2}\cdots
\xrightarrow{\F\alpha_{d-1}}\F A_d\xrightarrow{\F\alpha_d}\F A_{d+1}\xrightarrow{\zeta_{A_0}\F\alpha_{d+1}}\Sigma^d\F A_{0}$$
is also a $(d+2)$-angle since the first row is a $(d+2)$-angle. Then we are done.

In fact, since $\beta_1$ is a weak cokernel of $\F\alpha_1$ and $\F\alpha_1\circ\F\alpha_0=0$,
there exists a morphism $\omega_2\colon B_2\to\F A_2$ such that
$\F\alpha_1=\omega_2\beta_1$.
Since $\beta_2$ is a weak cokernel of $\beta_1$ and $\F\alpha_2\circ\omega_2\beta_1=\F\alpha_2\circ\F\alpha_1=0$,
there exists a morphism $\omega_3\colon B_3\to\F A_3$ such that
$\F\alpha_2\circ\omega_2=\omega_3\beta_2$.
Continue this process, there are morphisms $\omega_i\colon B_i\to \F A_i,~i=4,5,\cdots,d$
such that $\F\alpha_{i-1}\circ\omega_{i-1}=\omega_{i}\beta_{i-1}$.

It remains to show that there exists a morphism
$\omega_{d+1}\colon B_{d+1}\to\F A_{d+1}$
such that
\begin{equation}\label{eq7}
\begin{array}{l}
\zeta_{A_0}\F\alpha_{d+1}\circ\omega_{d+1}=\beta_{d+1}
\end{array}
\end{equation}
\begin{equation}\label{eq8}
\begin{array}{l}
\omega_{d+1}\beta_d=\F\alpha_d\circ\omega_d
\end{array}
\end{equation}

\textbf{Claim I}: The equality (\ref{eq7}) holds if and only if
\begin{equation}\label{eq9}
\begin{array}{l}
\eta_{A_{d+1}}(\omega_{d+1}\circ\Sigma^d\varphi\circ\alpha_{d+1})=(-1)^d
\eta_{A_0}(\Sigma^{-d}\beta_{d+1}\circ\varphi),~\mbox{for all}~ \varphi\in\Hom(A_0,\Sigma^{-d}B_{d+1}).
\end{array}
\end{equation}
The equality (\ref{eq8}) holds if and only if
\begin{equation}\label{eq10}
\begin{array}{l}
\eta_{A_{d+1}}(\omega_{d+1}\beta_d\phi)=
\eta_{A_d}(\omega_d\phi\alpha_d),~\mbox{for all}~\phi\in\Hom(A_{d+1},B_d).
\end{array}
\end{equation}
We first show that the equality (\ref{eq7}) holds if and only if the equality (\ref{eq9}) holds.

Suppose that the equality (\ref{eq7}) holds. For any morphism $\varphi\in\Hom(A_0,\Sigma^{-d}B_{d+1})$,
by the equality (\ref{eq7}), we have
$$\Sigma^{-d}(\zeta_{A_0}\F\alpha_{d+1}\circ\omega_{d+1})\circ\varphi=\Sigma^{-d}\beta_{d+1}\circ\varphi\colon A_0\to\F A_0.$$
It follows that
$$\begin{aligned}%
(-1)^d\eta_{A_0}(\Sigma^{-d}\beta_{d+1}\circ\varphi)&=(-1)^d\eta_{A_0}(\Sigma^{-d}(\zeta_{A_0}\F\alpha_{d+1}\circ\omega_{d+1})\circ\varphi)\\
                     &\xlongequal{(\ref{t2})}\eta_{\Sigma^dA_0}(\F\alpha_{d+1}\circ\omega_{d+1}\circ\Sigma^d\varphi)\\
                     &\xlongequal{(\ref{t1})}\eta_{A_{d+1}}(\omega_{d+1}\Sigma^d\varphi\circ\alpha_{d+1}).
\end{aligned}$$
This shows that the equality (\ref{eq9}) holds.

Conversely, assume that the equality (\ref{eq9}) holds.
For any morphism $\varphi\in\Hom(A_0,\Sigma^{-d}B_{d+1})$,
by the equalities (\ref{eq9}), (\ref{t1}) and (\ref{t2}), we have
$$\eta_{A_0}(\Sigma^{-d}\beta_{d+1}\circ\varphi)=\eta_{A_0}(\Sigma^{-d}(\zeta_{A_0}\F\alpha_{d+1}\circ\omega_{d+1})\circ\varphi).$$
 As  $$\Hom(A_0,\Sigma^{-d}B_{d+1})\times \Hom(\Sigma^{-d}B_{d+1}),\F A_0)\longrightarrow k,~~(\varphi,g)\longmapsto\eta_{A_0}(g\varphi)$$
is bilinear and non-degenerate, we have $\zeta_{A_0}\F\alpha_{d+1}\circ\omega_{d+1}=\beta_{d+1}$.
Namely, the equality (\ref{eq7}) holds.

Now we show that the equality (\ref{eq8}) holds if and only if the equality (\ref{eq10}) holds.

Assume that the equality (\ref{eq8}) holds. For any morphism $\phi\in\Hom(A_{d+1},B_d)$,
by the equality (\ref{eq8}), we have
$$\omega_{d+1}\beta_{d}\phi=\F\alpha_{d}\circ\omega_d\phi\colon A_{d+1}\to \F A_{d+1}.$$
It follows that
$$\eta_{A_{d+1}}(\omega_{d+1}\beta_{d}\phi)
=\eta_{A_{d+1}}(\F\alpha_{d}\circ\omega_d\phi)
\xlongequal{(\ref{t1})}\eta_{A_d}(\omega_d\phi\alpha_d).$$
This shows that the equality (\ref{eq10}) holds.

Conversely, assume that the equality (\ref{eq10}) holds.
For any morphism $\phi\in\Hom(A_{d+1},B_d)$, we have
$$\eta_{A_{d+1}}(\omega_{d+1}\beta_{d}\phi)\xlongequal{(\ref{eq10})}\eta_{A_d}(\omega_d\phi\alpha_d)
\xlongequal{(\ref{t1})}\eta_{A_{d+1}}(\F\alpha_{d}\circ\omega_d\phi).$$
By the bilinear map $$\Hom(A_{d+1},B_d)\times \Hom(B_d,\F A_{d+1})\longrightarrow k,~~(\phi,h)\longmapsto\eta_{A_{d+1}}(h\phi)$$
is non-degenerate, we obtain $\omega_{d+1}\beta_n=\F\alpha_d\circ\omega_d$.
That is to say, the equality (\ref{eq8}) holds.
\medskip

 \textbf{Claim II:} The following statements are equivalent.

$(\clubsuit)$~ The equalities (\ref{eq9}) and (\ref{eq10}) hold.

 $(\spadesuit)$~ If there are morphisms $\varphi\colon A_0\to \Sigma^{-d}B_{d+1}$ and $\phi\colon A_{d+1}\to B_d$ such that
$\Sigma^d\varphi\circ\alpha_{d+1}=\beta_d\phi$, then
\begin{equation}\label{eq11}
\begin{array}{l}
\eta_{A_{d}}(\omega_{d}\phi\alpha_d)=(-1)^d\eta_{A_0}(\Sigma^{-d}\beta_{d+1}\circ\varphi).
\end{array}
\end{equation}
If $(\clubsuit)$  hold, it is obvious that the equality (\ref{eq11}) also holds.

Conversely, assume that $(\spadesuit)$  holds. Put
$$V_1:=\{\Sigma^d\varphi\circ\alpha_{d+1}~|~\varphi\in\Hom(A_0,\Sigma^{-d}B_{d+1})\}\subseteq\Hom(A_{d+1},B_{d+1}),$$
$$V_2:=\{\beta_d\phi~|~\phi\in\Hom(A_{d+1},B_d)\}\subseteq\Hom(A_{d+1},B_{d+1}).$$
Take a set of basis elements $a_1,a_2,\cdots,a_m$ of $V_1\cap V_2$, and expand them to a set of basis elements
$$a_1,\cdots,a_m,a_{m+1},\cdots,a_{m+t}$$
 of $V_1$. At the same time, it is also extended to a set of basis elements
 $$b_1=a_1,\cdots,b_m=a_m,b_{m+1},\cdots,b_{m+s}$$
 of $V_2$.
Take a set of basis elements $c_1,\cdots,c_n$ of $\Hom(B_{n+1},\F A_{n+1})$.
Write
$$\eta_{A_{n+1}}(c_ja_i)=a_{ij},~i=1,\cdots,m+t;~j=1,\cdots,n,$$
$$\eta_{A_{n+1}}(c_jbi)=b_{ij},~i=1,\cdots,m+s;~j=1,\cdots,n,$$
$$M=(a_{ij})_{(m+t)\times n},~N=(b_{ij})_{(m+s)\times n}.$$
 Thus the first $m$ rows of the matrices of $M$ and $N$ are equal.
Since $a_1,\cdots,a_m,a_{m+1},\cdots,a_{m+t}$ is a set of  basis elements of $V_1$, there are
morphisms $$\varphi_1,\cdots,\varphi_m,\varphi_{m+1},\cdots,\varphi_{m+t}\in\Hom(A_0,\Sigma^{-d}B_{d+1})$$
such that
$$a_1=\Sigma^{d}\varphi_1\circ\alpha_{d+1},\cdots,a_m=\Sigma^{d}\varphi_m\circ\alpha_{d+1},
a_{m+1}=\Sigma^{d}\varphi_{m+1}\circ\alpha_{d+1},\cdots,
a_{m+t}=\Sigma^{d}\varphi_{m+t}\circ\alpha_{d+1}.$$
Put
$$d_i=(-1)^d\eta_{A_0}(\Sigma^{-d}\beta_{d+1}\circ\varphi),~~i=1,\cdots,m+t.$$
Similarly, there are morphisms
$$\phi_1,\cdots,\phi_m,\phi_{m+1},\cdots,\phi_{m+t}$$
such that
$$b_1=\beta_d\phi_1,\cdots,b_m=\beta_d\phi_m,\cdots,b_{m+1}=\beta_d\phi_{m+1},\cdots,b_{m+s}=\beta_d\phi_{m+s}.$$
Put
$$e_i=\eta_{A_n}(\omega_d\phi_i\alpha_d),~~i=1,\cdots,m+s.$$
Since $\beta_d\phi_i=b_i=a_i=\Sigma^d\varphi_i\circ\alpha_{d+1},~~i=1,\cdots,m.$
By the equality (\ref{eq11}), we have
$$e_i=\eta_{A_d}(\omega_d\phi_i\alpha_d)=(-1)^d(\Sigma^{-d}\beta_{d+1}\circ\varphi_i)=d_i,~~i=1,\cdots,m.$$
It follows that there exists a morphism
$\omega_{d+1}\colon B_{d+1}\to \F A_d$ such that the equalities (\ref{eq9}) and (\ref{eq10}) hold
if and only if the following two linear equations have solutions
$$M\left(\begin{smallmatrix}
x_1\\\vdots\\x_n
              \end{smallmatrix}
            \right)=\left(\begin{smallmatrix}
d_1\\\vdots\\d_{m+t}
              \end{smallmatrix}
            \right),$$
$$N\left(\begin{smallmatrix}
x_1\\\vdots\\x_n
              \end{smallmatrix}
            \right)=\left(\begin{smallmatrix}
e_1\\\vdots\\e_{m+s}
              \end{smallmatrix}
            \right).$$
 Note that the first $m$ rows of the matrices $M$ and $N$ are equal.
Let $L=\left(\begin{smallmatrix}
M\\N'\end{smallmatrix}\right)$
where $N'$ consists of the last $s$ rows of $N$.
Hence, the necessary and sufficient  conditions for the two linear equations
 stated above  to have a solution is  that the following linear equations
\begin{equation}\label{eq12}
\begin{array}{l}
L\left(\begin{smallmatrix}
x_1\\\vdots\\x_n
              \end{smallmatrix}
            \right)=\left(\begin{smallmatrix}
d_1\\\vdots\\d_m\\d_{m+1}\\ \vdots\\
d_{m+t}\\
e_{m+1}\\
\vdots\\
e_{m+s}
              \end{smallmatrix}
            \right).\end{array}
\end{equation}
have a solution. Since the bilinear map $$\Hom(A_{d+1},B_d)\times \Hom(B_d,\F A_{d+1})\longrightarrow k,~~(u,v)\longmapsto\eta_{A_{d+1}}(vu)$$
is non-degenerate, we get that $L$ has row full rank, and thus, that  (\ref{eq12}) has a solution.

\textbf{Step 3:}  We show that the equality (\ref{eq11}) holds.

Assume that there are morphisms $\varphi\colon A_0\to \Sigma^{-d}B_{d+1}$ and $\phi\colon A_{d+1}\to B_d$ such that $\Sigma^d\varphi\circ\alpha_{d+1}=\beta_d\phi.$
Consider the following morphisms of $(d+2)$-angles in $\C$:
$$\xymatrix@C=1.2cm{A_1\ar[r]^{\alpha_1}\ar@{-->}[d]^{\delta_1}&A_2\ar[r]^{\alpha_2}\ar@{-->}[d]^{\delta_2}&A_3\ar[r]^{\alpha_3}
\ar@{-->}[d]^{\delta_3}
&\cdot\cdot\cdot\ar[r]^{\alpha_{d-1}}
&A_d\ar[r]^{\alpha_d}\ar@{-->}[d]^{\delta_d}&A_{d+1}\ar[r]^{\alpha_{d+1}}
\ar[d]^{\phi}&\Sigma^d A_0\ar[d]^{\Sigma^d\varphi}\ar[r]^{(-1)^d\Sigma^d\alpha_0}&\Sigma^d A_1\ar@{-->}[d]^{\Sigma^d\delta_1}\\
\F A_0\ar[r]^{\F\alpha_0}&\F A_1\ar[r]^{\beta_1}&B_2\ar[r]^{\beta_2}&\cdot\cdot\cdot\ar[r]^{\beta_{d-2}}&B_{d-1}\ar[r]^{\beta_{d-1}}& B_d\ar[r]^{\beta_d}&B_{d+1}\ar[r]^{\beta_{d+1}}&\Sigma^d\F A_{0}}$$
where $\delta_1,\delta_2,\cdots,\delta_d$ exist by the axiom (N3). It follows that
$$\begin{aligned}%
\eta_{A_{d}}(\omega_{d}\phi\alpha_d)&=\eta_{A_{d}}(\omega_{d}\beta_{d-1}\delta_d)=\eta_{A_{d}}(\F\alpha_{d-1}\circ\omega_{d-1}\delta_d)\\
                     &\xlongequal{(\ref{t1})}\eta_{A_{d-1}}(\omega_{d-1}\delta_d\alpha_{d-1})=\eta_{A_{d-1}}(\omega_{d-1}\beta_{d-2}\delta_{d-1})
                     =\eta_{A_{d-1}}(\F\alpha_{d-2}\circ\omega_{d-2}\delta_{d-1})\\
                     &\xlongequal{(\ref{t1})}\eta_{A_{d-2}}(\omega_{d-2}\delta_{d-1}\alpha_{d-2})\\
                     &\vdots\\
&\xlongequal{(\ref{t1})}\eta_{A_{2}}(\omega_{2}\delta_{3}\alpha_{2})=\eta_{A_{2}}(\omega_{2}\beta_{1}\delta_{2})=\eta_{A_{2}}(\F\alpha_{1}\circ\delta_{2})\\
&\xlongequal{(\ref{t1})}\eta_{A_{1}}(\delta_{2}\alpha_{1})=\eta_{A_{1}}(\F\alpha_0\circ\delta_{1})\\
&\xlongequal{(\ref{t1})}\eta_{A_{0}}(\delta_{1}\alpha_{0})=(-1)^d\eta_{A_0}(\Sigma^{-d}\beta_{d+1}\circ\varphi).
\end{aligned}$$
as we wished. This completes the proof.  \qed

\section{Connection between Auslander-Reiten $(d+2)$-angles and Serre functors}
\setcounter{equation}{0}
In this section, we build a link between Auslander-Reiten $(d+2)$-angles and Serre functors.

\begin{lemma}\label{lem1}
Let $\C$ be a  $(d+2)$-angulated category,
$$A_{\bullet}:~~\xymatrix {A_0 \xrightarrow{~\alpha_0~}A_1 \xrightarrow{~\alpha_1~} A_2 \xrightarrow{~\alpha_2~} \cdots
  \xrightarrow{~\alpha_{d - 1}~} A_d \xrightarrow{~\alpha_{d}~} A_{d+1}\xrightarrow{~\alpha_{d+1}~} \Sigma^d A_0}$$
and
$$B_{\bullet}:~~\xymatrix {B_0 \xrightarrow{~\beta_0~}B_1 \xrightarrow{~\beta_1~} B_2 \xrightarrow{~\beta_2~} \cdots\xrightarrow{~\beta_{d - 1}~} B_d \xrightarrow{~\beta_{d}~} B_{d+1}\xrightarrow{~\beta_{d+1}~} \Sigma^d B_0}$$
be two Auslander-Reiten $(d+2)$-angles in $\C$.
Then $A_0\simeq B_0$ if and only if $A_{d+1}\simeq B_{d+1}$.
\end{lemma}

\proof Assume that $\varphi_{d+1}\colon A_{d+1}\to B_{d+1}$ is an isomorphism.
Then $\varphi_{d+1}\alpha_d$ is not split epimorphism, otherwise, $\alpha_d$
is a split epimorphism and this is impossible. Since $\alpha_d$ is right almost split, there exists  a morphism $\varphi_d\colon A_d\to B_d$ such that
$\varphi_{d+1}\alpha_d=\beta_d\varphi_d$. Thus we have the following commutative diagram
$$\xymatrix@C=1.2cm{
A_0\ar[r]^{\alpha_0}\ar@{-->}[d]^{\varphi_0}&A_1 \ar[r]^{\alpha_1}\ar@{-->}[d]^{\varphi_1}  & A_2 \ar[r]^{\alpha_2}\ar@{-->}[d]^{\varphi_2} & \cdots \ar[r]^{\alpha_{d-1}}& A_d \ar[r]^{\alpha_d\;\;}\ar[d]^{\varphi_{d}}&A_{d+1}\ar[r]^{\alpha_{d+1}}\ar[d]^{\varphi_{d+1}\;} & \Sigma^dA_0 \ar@{-->}[d]^{\Sigma^d\varphi_0}\\
B_0\ar[r]^{\beta_0}&B_1 \ar[r]^{\beta_1} & B_2\ar[r]^{\beta_2} & \cdots \ar[r]^{\beta_{d-1}} & B_d \ar[r]^{\beta_d\;\;}&B_{d+1}\ar[r]^{\beta_{d+1}\;}& \Sigma^dB_0\\
}$$
of $(d+2)$-angles in $\C$.

 Since $A_{\bullet}$ and $B_{\bullet}$ are Auslander-Reiten $(d+2)$-angles, by Remark \ref{rem},
we have that $A_0$ and $B_0$ are indecomposable objects.  If $\varphi_0$ is not an isomorphism,
we obtain that $\varphi_0$ is not a split monomorphism since $A_0$ and $B_0$ are indecomposable. Because $\alpha_0$ is left
almost split, there exists a morphism $h\colon A_1\to B_0$ such that
$\varphi_0=h\alpha_0$.
Hence $\beta_{d+1}\varphi_{d+1}=(\Sigma^d\varphi_0)\alpha_{d+1}=\Sigma^dh(\Sigma^d\alpha_0)\alpha_{d+1}=0$
which implies  $\beta_{d+1}=0$. This is a contradiction since $B_{\bullet}$ is an
Auslander-Reiten $(d+2)$-angles in $\C$.
Hence $\varphi_0\colon A_0\to B_0$ is an isomorphism.  \qed

\begin{remark}
Let $\C$ be a $(d+2)$-angulated category, and
$$\xymatrix {A_0 \xrightarrow{~\alpha_0~}A_1 \xrightarrow{~\alpha_1~} A_2 \xrightarrow{~\alpha_2~} \cdots
  \xrightarrow{~\alpha_{d - 1}~} A_d \xrightarrow{~\alpha_{d}~} A_{d+1}\xrightarrow{~\alpha_{d+1}~} \Sigma^d A_0}$$
be an Auslander-Reiten $(d+2)$-angle in $\C$.
By Lemma \ref{lem1}, we know that $A_0$ is uniquely determined up to isomorphism.
In this case, we  write $A_0=\tau_dA_{d+1}$.
\end{remark}
\medskip

In order to prove the main result of this section, we need the following lemma.

\begin{lemma}\label{lem3}
Let $\C$ be a $(d+2)$-angulated category. Suppose that
$$\xymatrix {\tau_dC \xrightarrow{~\alpha_0~}C_1 \xrightarrow{~\alpha_1~} C_2 \xrightarrow{~\alpha_2~} \cdots
  \xrightarrow{~\alpha_{d - 1}~} C_d \xrightarrow{~\alpha_{d}~} C\xrightarrow{~\alpha_{d+1}~} \Sigma^d\tau_dC}$$
is an Auslander-Reiten $(d+2)$-angle in $\C$ and $B$ is a indecomposable object in $\C$. Then the following hold:
\begin{itemize}
\item[\emph{(1)}] For any non-zero $g\in\emph{\Hom}_{\C}(B,\tau_d\Sigma^dC)$, there is an $f\in\emph{Hom}_{\C}(C,B)$ such that $\alpha_{d+1}=gf$.

\item[\emph{(2)}] For any non-zero $f\in\emph{\Hom}_{\C}(C,B)$, there is a $g\in\emph{\Hom}_{\C}(B,\tau_d\Sigma^dC)$ such that $\alpha_{d+1}=gf$.
\end{itemize}
\end{lemma}

\proof (1) For any non-zero morphism $g\colon B\to\tau_d\Sigma^dC$, there exists
a $(d+2)$-angle
$$\xymatrix {\tau_dC \xrightarrow{~\beta_0~}B_1 \xrightarrow{~\beta_1~} B_2 \xrightarrow{~\beta_2~} \cdots
  \xrightarrow{~\beta_{d - 1}~} B_d \xrightarrow{~\beta_{d}~} B\xrightarrow{~g~} \Sigma^d\tau_dC}$$
in $\C$. Since $g$ is non-zero, then $\beta_0$ is not split monomorphism. It follows that there exists a morphism $h\colon C_1\to B_1$ such that $\beta_0=h\alpha_0$.
Thus we have the following commutative diagram
$$\xymatrix@C=1.2cm{
\tau_dC\ar[r]^{\alpha_0}\ar@{=}[d]&C_1 \ar[r]^{\alpha_1}\ar[d]^{h}  & C_2 \ar[r]^{\alpha_2}\ar@{-->}[d]^{\varphi_2} & \cdots \ar[r]^{\alpha_{d-1}}& C_d \ar[r]^{\alpha_d\;\;}\ar@{-->}[d]^{\varphi_{d}}&C\ar[r]^{\alpha_{d+1}\quad}\ar@{-->}[d]^{f} & \Sigma^d\tau_dC \ar@{=}[d]\\
\tau_dC\ar[r]^{\beta_0}&B_1 \ar[r]^{\beta_1} & B_2\ar[r]^{\beta_2} & \cdots \ar[r]^{\beta_{d-1}} & B_d \ar[r]^{\beta_d\;\;}&B\ar[r]^{g\quad}& \Sigma^d\tau_dC
}$$
of $(d+2)$-angles in $\C$. We obtain $\alpha_{d+1}=gf$.

(2) For any non-zero morphism $f\colon C\to B$, there exists a
$(d+2)$-angle
$$\xymatrix {A_0 \xrightarrow{~\gamma_0~}A_1 \xrightarrow{~\gamma_1~} \cdots
  \xrightarrow{~\gamma_{d-2}~}A_{d-1}\xrightarrow{~\gamma_{d-1}~} C \xrightarrow{~f~} B\xrightarrow{~\gamma~} \Sigma^d A_0}$$
in $\C$. Since $f\gamma_{d-1}=0$ and $f\neq 0$, we have that $\gamma_{d-1}$ is not
split epimorphism.
Since $\alpha_d$ is right almost split, there exists a morphism
$u\colon A_{d-1}\to C_d$ such that $\alpha_du=\gamma_{d-1}$
 and then $$\alpha_{d+1}\gamma_{d-1}=\alpha_{d+1}\alpha_du=0.$$
So there exists a morphism $g\colon B\to \Sigma^d\tau_dC$ such that
$\alpha_{d+1}=gf$.  \qed
\medskip

The following lemma can be found in \cite[Proposition I.1.4]{rv}.

\begin{lemma}\label{lem2}
Let $\C$ be an additive category and $\F\colon \C\to\C$ be an additive functor.
If there exists a $k$-linear map
$\eta_X\colon{\rm Hom}(X,\F X)\longrightarrow k$ such that
the bilinear map
$$(-,-)\colon{\rm Hom}(X,Y)\times{\rm Hom}(Y,\F X)\longrightarrow k,~~(f,g)=\eta_X(gf)$$
is a non-degenerate pairing for all $X,Y\in\C$, then
$\F$ can be viewed as a right Serre functor, where
$$\eta_{X,Y}\colon{\rm Hom}(X,Y)\longrightarrow {\rm Hom}(Y,\F Y)^{\ast}$$
given by
$$\eta_{X,Y}(f)(g)=\eta_X(gf),~\mbox{for all}~ f \in\emph{\Hom}(X,Y), g\in\emph{\Hom}(Y,\F X).$$
\end{lemma}

\medskip
Now we state and prove our second main result. This generalizes the work by Reiten and Van den Bergh \cite[Theorem I.2.4]{rv} for triangulated categories.

\begin{theorem}\label{main2}
Let $\C$ be a $(d+2)$-angulated category with  $d$-suspension
functor $\Sigma^d$.
Then $\C$ has Auslander-Reiten $(d+2)$-angles if and only if
$\C$ has a Serre functor $\mathbb{S}$.

If either of these properties holds, then the action of the Serre functor on objects coincides with $\tau_d\Sigma^d$, namely
$\tau_d=\mathbb{S}\Sigma^{-d}$. In this case, $\tau_d$ is called  the $d$-Auslander-Reiten translation.
\end{theorem}

\proof We first show the `only if' part. Suppose that $\C$ has  Auslander-Reiten $(d+2)$-angles. For any indecomposable object $C\in\C$, put $\F C=\tau_d \Sigma^d C$. Then there exists an Auslander-Reiten $(d+2)$-angle
$$\xymatrix {\tau_dC \xrightarrow{~\alpha_0~}C_1 \xrightarrow{~\alpha_1~} C_2 \xrightarrow{~\alpha_2~} \cdots
  \xrightarrow{~\alpha_{d - 1}~} C_d \xrightarrow{~\alpha_{d}~} C\xrightarrow{~\alpha_{d+1}~} \Sigma^d\tau_dC.}$$
Thus $\alpha_{d+1}\neq 0$. For the indecomposable object $C$, take a $k$-linear map
$$\eta_{C}\colon \Hom(C,\F C)\to k$$
such that $\eta_{C}(\alpha_{d+1})\neq0$.
By Lemma \ref{lem3}, for any indecomposable object $B\in\C$, the pairing
$$\Hom(C,B)\times\Hom(B,\F C)\to k \ \textrm{given by}\ (f,g)=\eta_{C}(gf)$$
is a non-degenerate $k$-bilinear map.
 By Lemma \ref{lem2}, $\F$ can be viewed as a right Serre functor.
Dually, we can prove that $\C$ has a left Serre functor.

We now prove the other implication.
 Assume that $\C$ has a Serre functor $\mathbb{S}$. Then for any indecomposable object
$C$, there exists a $k$-linear isomorphism
$$(\eta^{\ast}_{C,C})^{-1}\colon \End(C)^{\ast}\to \Hom(C,\F C).$$
Since $\End(C)$ is local and $k$ is an algebraically closed field, we have
$\End(C)/\rad\End(C)\simeq k$.
 Then $\theta_C\colon \End(C)\to \End(C)/\rad\End(C)\simeq k$ is a natural projection  with $\Ker~\theta_C=\rad\End(C)$.
For $$\delta=(\eta^{\ast}_{C,C})^{-1}(\theta_C)\in\Hom(C,\F C),$$
we have $\delta\neq0$ and  so $\alpha_d$ is not a split epimorphism. Consider a $(d+2)$-angle extending $\delta$:
$$\xymatrix {\Sigma^{-d}\F C \xrightarrow{~\alpha_0~}C_1 \xrightarrow{~\alpha_1~} C_2 \xrightarrow{~\alpha_2~} \cdots
  \xrightarrow{~\alpha_{d-2}~} C_{d-1}\xrightarrow{~\alpha_{d-1}~} C_d \xrightarrow{~\alpha_d~} C\xrightarrow{~\delta~} \F C}$$
where when $d\geq2$, we can choose $\alpha_1,\alpha_2,\cdots,\alpha_{d-1}$ in
$\rad_{\C}$ by Lemma \ref{lemma}.

We claim that
$$\xymatrix {\Sigma^{-d}\F C \xrightarrow{~\alpha_0~}C_1 \xrightarrow{~\alpha_1~} C_2 \xrightarrow{~\alpha_2~} \cdots
  \xrightarrow{~\alpha_{d-2}~} C_{d-1}\xrightarrow{~\alpha_{d-1}~} C_d \xrightarrow{~\alpha_d~} C\xrightarrow{~\delta~} \F C}$$
is  an Auslander-Reiten $(d+2)$-angle in $\C$.
 We first show that $\alpha_d$ is a right almost split morphism.
Let $\beta\colon B\to C$  not be a split epimorphism in $\C$.
Then the composition
$$\Hom(C,B)\xrightarrow{\Hom(C,\ \beta)}\End(C)\xrightarrow{~\theta_Z~}k$$
is zero. Consider the following commutative diagram
$$\xymatrix@C=1.2cm {\Hom(C,C)^\ast\ar[r]^{(\eta^{\ast}_{C,C})^{-1}}\ar[d]_{\Hom(C,\ \beta )^\ast}&\Hom(C,\F C)\ar[d]^{\Hom(\beta, \ \F C)}\\
\Hom(C,B)^\ast\ar[r]^{(\eta^{\ast}_{C,B})^{-1}}& \Hom(B,\F C)}$$
we can obtain
$$\delta\beta=(\eta^{\ast}_{C,C})^{-1}(\theta_Z)\beta=(\eta^{\ast}_{C,B})^{-1}\Hom(C, \beta )^\ast(\theta_C)=(\eta^{\ast}_{C,B})^{-1}(\theta_C)\Hom(C, \beta)=0.$$
So there exists a morphism $h\colon C\to C_d$ such that $\beta=\alpha_dh$.
This shows that $\alpha_d$ is a right almost split morphism.

Since $C$ is an indecomposable object, we have that
$\Sigma^{-d}\F C$ is also indecomposable which implies that $\End(\Sigma^{-d}\F C)$
is local. By  Lemma \ref{keylemma}, we know that
$$\xymatrix {\Sigma^{-d}\F C \xrightarrow{~\alpha_0~}C_1 \xrightarrow{~\alpha_1~} C_2 \xrightarrow{~\alpha_2~} \cdots
  \xrightarrow{~\alpha_{d-2}~} C_{d-1}\xrightarrow{~\alpha_{d-1}~} C_d \xrightarrow{~\alpha_d~} C\xrightarrow{~\delta~} \F C}$$
is an Auslander-Reiten $(d+2)$-angle in $\C$.

Dually, we can show that for any indecomposable object $C\in\C$, there exists
an Auslander-Reiten $(d+2)$-angle
$\xymatrix {C \xrightarrow{~\beta_0~}C_1 \xrightarrow{~\beta_1~} C_2 \xrightarrow{~\beta_2~} \cdots
  \xrightarrow{~\beta_{d-1}~} C_d \xrightarrow{~\beta_d~} \mathbb{S}^{-1}\Sigma^dC\xrightarrow{~\beta_{d+1}~} \Sigma^dC}$.  \qed
\medskip

Now we give an example illustrating our main result in this section.

\begin{example}
We first recall the standard construction of $(d+2)$-angulated categories given by Geiss-Keller-Oppermann \cite[Theorem 1]{gko}.
Let $\C$ be a triangulated category and $\mathcal{T}$ a $d$-cluster tilting subcategory which is closed under $\Sigma^{d}$, where $\Sigma$ is the shift functor of $\C$. Then $(\mathcal{T},\Sigma^{d},\Theta)$ is a $(d+2)$-angulated category, where $\Theta$ is the class of all sequences
$$A_0\xrightarrow{f_0}A_1\xrightarrow{f_1}A_2\xrightarrow{f_2}\cdots\xrightarrow{f_{d-1}}A_d\xrightarrow{f_d}A_{d+1}\xrightarrow{f_{d+1}}\Sigma^{d} A_0$$
such that there exists a diagram
$$\xymatrixcolsep{0.4pc}
 \xymatrix{& A_1 \ar[dr]\ar[rr]^{f_1}  &  & A_2  \ar[dr]  & & \cdots  & & A_{d} \ar[dr]^{f_{d}}      \\
A_0 \ar[ur]^{f_0} & \mid & \ar[ll]  A_{1.5}\ar[ur] & \mid &  \ar[ll]  A_{2.5} & \cdots & A_{d-1.5}\ar[ur] & \mid & \ar[ll] A_{d+1}   }$$
with $A_i\in\mathcal{T}$ for all $i\in\mathbb{Z}$, such that all oriented triangles are triangles in $\C$, all non-oriented triangles commute, and $f_{d+1}$ is the composition along the lower edge of the diagram.

Assume that $\C$ has a Serre functor $\mathbb{S}$.  By \cite[Proposition 3.4]{iy},
we know that $\mathcal{T}=\mathbb{S}\Sigma^{-d}\mathcal T=\Sigma^d\mathbb{S}^{-1}\mathcal T$.
By hypothesis, $\Sigma^d\mathcal T=\mathcal T$.
It follows that $\mathbb{S}\mathcal T=\mathcal T=\mathbb{S}^{-1}\mathcal T$.
 In particular, $\mathcal T$ has an auto-equivalence functor $\mathbb{S}$.
Hence $(\mathcal{T},\Sigma^{d},\Theta)$ has
a Serre functor $\mathbb{S}$. By Theorem \ref{main2}, we know that  $\mathcal T$ has Auslander-Reiten $(d+2)$-angles and $\tau_d=\mathbb{S}\Sigma^{-d}$.
\end{example}

\section{$(d+2)$-angulated quotient categories}
Let $\C$ be  an additive category,
and $\X$ be a subcategory of $\C$.
Recall that we say a morphism $f\colon A \to B$ in $\C$ is an $\X$-\emph{monic} if
$$\Hom_{\C}(f,X)\colon \Hom_{\C}(B,X) \to \Hom_{\C}(A,X)$$
is an epimorphism for all $X\in\X$. We say that $f$ is an $\X$-\emph{epic} if
$$\Hom_{\C}(X,f)\colon \Hom_{\C}(X,A) \to \Hom_{\C}(X,B)$$
is an epimorphism for all $X\in\X$.
 Furthermore,
we say that $f$ is a left $\X$-approximation of $A$ if $f$ is an $\X$-monoic and $B\in\X$.
We say that $f$ is a right $\X$-approximation of $B$ if $f$ is an $\X$-epic and $A\in\X$.

A subcategory $\X$ is called \emph{contravariantly finite} if any object in $\C$ admits a right
$\X$-approximation. Dually we can define  \emph{covariantly finite} subcategory.

\begin{definition}\label{dd1}
Let $\C$ be a $(d+2)$-angulated category. A subcategory $\X$ of $\C$ is called
\emph{strongly contravariantly finite}, if for any object $C\in\C$, there exists a $(d+2)$-angle
$$B\xrightarrow{}X_1\xrightarrow{}X_2\xrightarrow{}\cdots\xrightarrow{}X_{d-1}\xrightarrow{}X_{d}\xrightarrow{~g~}C\xrightarrow{~~}\Sigma^dB$$
where $g$ is a right $\X$-approximation of $C$ and $X_1,\cdots,X_d\in\X$.
Dually we can define a \emph{strongly covariantly finite} subcategory.

A strongly contravariantly finite and strongly  covariantly finite subcategory is called \emph{ strongly functorially finite}.
\end{definition}
\medskip

\begin{lemma}\label{lem0}
Let $\C$ be a $(d+2)$-angulated category and $\X$ be a strongly covariantly finite subcategory of $\C$.  Then the quotient
category $\C/\X$ is a right $(d+2)$-angulated category with the following endofunctor $\mathbb{H}$ and right $(d+2)$-angles:
\begin{itemize}
\item[\emph{(1)}] For any object $A_0\in\C$, we take a $(d+2)$-angle
$$A_0\xrightarrow{~f_0~}X_1\xrightarrow{~f_1~}X_2\xrightarrow{~f_2~}\cdots\xrightarrow{~f_{d-1}~}X_d\xrightarrow{~f_d~} \mathbb{H}A_0\xrightarrow{~f_{d+1}~}\Sigma^d A_0$$
with $f_0$ is a left $\X$-approximation of $A_0$ and $X_1,X_2,\cdots,X_d\in\X$. Then $\mathbb{H}$ gives a well-defined endofunctor of $\C/\X$.
\item[\emph{(2)}] For any $(d+2)$-angle
 $$A_0\xrightarrow{~g_0~}A_1\xrightarrow{~g_1~}A_2\xrightarrow{~g_2~}A_3\xrightarrow{~g_3~}\cdots\xrightarrow{~g_{d-1~}}A_d\xrightarrow{~g_d~}A_{d+1}\xrightarrow{~g_{d+1}~}\Sigma^d A_0$$
with $g_0$ is an $\X$-monic, take the following commutative diagram of $(d+2)$-angles
$$\xymatrix{
A_0 \ar[r]^{g_0}\ar@{=}[d]& A_1 \ar[r]^{g_1}\ar@{-->}[d]^{\varphi_1} & A_2 \ar[r]^{g_2}\ar@{-->}[d]^{\varphi_2}  & \cdots \ar[r]^{g_{d-1}}& A_d \ar[r]^{g_d}\ar@{-->}[d]^{\varphi_d}&A_{d+1}\ar[r]^{g_{d+1}}\ar@{-->}[d]^{\varphi_{d+1}} & \Sigma^d A_0 \ar@{=}[d]\\
A_0 \ar[r]^{f_0}&X_1 \ar[r]^{f_1} & X_2 \ar[r]^{f_2}  & \cdots \ar[r]^{f_{d-1}} & X_d \ar[r]^{f_d}&\mathbb{H}A_0\ar[r]^{f_{d+1}}& \Sigma^d A_0.\\
}$$
Then we have a complex
$$A_0\xrightarrow{~\overline{g_0}~} A_1\xrightarrow{~\overline{g_1}~}A_2\xrightarrow
{~\overline{g_2}~}\cdots\xrightarrow{~\overline{g_{d-1}}~}A_{d}\xrightarrow{~\overline{g_{d}}~}A_{d+1}\xrightarrow{~\overline{\varphi_{d+1}}~}\mathbb{H}A_0.$$ We define right $(d+2)$-angles in $\C/\X$ as the complexes which are isomorphic to complexes obtained in this way.
\end{itemize}
\end{lemma}

\proof
Since the proof is similar to \cite[Theorem 3.7]{l1}, we omit it. See also \cite[Remark 3.8]{l2}. \qed
\medskip

 The following is just the dual statement of Lemma \ref{lem0}.

\begin{lemma}\label{duallemma}
Let $\C$ be a $(d+2)$-angulated category and $\X$ be a strongly contravariantly finite subcategory of $\C$.  Then the quotient
category $\C/\X$ is a left $(d+2)$-angulated category with the following endofunctor $\mathbb{L}$ and left $(d+2)$-angles:
\begin{itemize}
\item[\emph{(1)}] For any object $A_{d+1}\in\C$, we take a $(d+2)$-angle
$$\mathbb{L}A_{d+1}\xrightarrow{~f_0~}X_1\xrightarrow{~f_1~}X_2\xrightarrow{~f_2~}\cdots\xrightarrow{~f_{d-1}~}X_d\xrightarrow{~f_d~} A_{d+1}\xrightarrow{~f_{d+1}~}\Sigma^d\mathbb{L}A_{d+1}$$
with $f_d$ is a right $\X$-approximation of $A_{d+1}$ and $X_1,X_2,\cdots,X_d\in\X$. Then $\mathbb{L}$ gives a well-defined endofunctor of $\C/\X$.
\item[\emph{(2)}] For any $(d+2)$-angle
 $$A_0\xrightarrow{~g_0~}A_1\xrightarrow{~g_1~}A_2\xrightarrow{~g_2~}A_3\xrightarrow{~g_3~}\cdots\xrightarrow{~g_{d-1~}}A_d\xrightarrow{~g_d~}A_{d+1}\xrightarrow{~g_{d+1}~}\Sigma^d A_0$$
with $g_d$ is an $\X$-epic, take the following commutative diagram of $(d+2)$-angles
$$\xymatrix{
\mathbb{L}A_{d+1} \ar[r]^{\;\;f_0}\ar@{-->}[d]^{\varphi_0}& X_1 \ar[r]^{f_1}\ar@{-->}[d]^{\varphi_1} & X_2 \ar[r]^{f_2}\ar@{-->}[d]^{\varphi_2}  & \cdots \ar[r]^{f_{d-1}}& X_d \ar[r]^{f_d}\ar@{-->}[d]^{\varphi_d}&A_{d+1}\ar[r]^{f_{d+1}\quad}\ar@{=}[d] & \Sigma^d \mathbb{L}A_{d+1} \ar@{-->}[d]^{\Sigma^d\varphi_0}\\
A_0 \ar[r]^{g_0}&A_1 \ar[r]^{g_1} & A_2 \ar[r]^{g_2}  & \cdots \ar[r]^{g_{d-1}} & A_d \ar[r]^{g_d}&A_{d+1}\ar[r]^{g_{d+1}}& \Sigma^d A_0.\\
}$$
Then we have a complex
$$\mathbb{L}A_{d+1}\xrightarrow{~\overline{\varphi_0}~} A_0\xrightarrow
{~\overline{g_0}~} A_1\xrightarrow{~\overline{g_1}~}A_2\xrightarrow
{~\overline{g_2}~}\cdots\xrightarrow{~\overline{g_{d-1}}~}A_{d}\xrightarrow{~\overline{g_{d}}~}A_{d+1}.$$ We define left $(d+2)$-angles in $\C/\X$ as the complexes which are isomorphic to complexes obtained in this way.
\end{itemize}
\end{lemma}

\medskip
The notion of mutation pairs of subcategories in a $(d+2)$-angulated category was
defined by Lin \cite[Definition 3.1]{l1}. We recall the definition here.

\begin{definition}
Let $\C$ be a $(d+2)$-angulated category, and $\X\subseteq\A$ be two subcategories of $\C$.
 The pair $(\A,\A)$ is called an $\X$-\emph{mutation pair} if it satisfies the following conditions:
\begin{itemize}
\item[(1)] For any object $A\in\A$, there exists a $(d+2)$-angle $$A\xrightarrow{x_0}X_1\xrightarrow{x_1}X_2\xrightarrow{x_2}\cdots\xrightarrow{x_{d-2}}X_{d-1}\xrightarrow{x_{d-1}}X_{d}\xrightarrow{x_{d}}B\xrightarrow{x_{d+1}}\Sigma^d A$$
where $X_i\in\X, B\in\A, x_0$ is a left $\X$-approximation of $A$ and $x_{d}$ is a right $\X$-approximation of $B$.

\item[(2)] For any object $C\in\A$, there exists a $(d+2)$-angle $$D\xrightarrow{x'_0}X'_1\xrightarrow{x'_1}X'_2\xrightarrow{x'_2}\cdots\xrightarrow{x'_{d-2}}X'_{d-1}\xrightarrow{x'_{d-1}}X'_d\xrightarrow{x'_{d}}C\xrightarrow{x'_{d+1}}\Sigma^d D$$
where $X'_i\in\X, D\in\A, x'_0$ is a left $\X$-approximation of $D$ and $x'_{d}$ is a right $\X$-approximation of $C$.
\end{itemize}
\end{definition}

Let $\C$ be a $(d+2)$-angulated category. Recall that a subcategory $\A$ of $\C$ is called
\emph{extension closed} if for any morphism $\alpha_{d+1}\colon A_{d+1}\to\Sigma A_0$ with
$A_0,A_{d+1}\in\A$, there exists a $(d+2)$-angle
$$A_0\xrightarrow{\alpha_0}A_1\xrightarrow{\alpha_1}A_2\xrightarrow{\alpha_2}\cdots\xrightarrow{\alpha_{d-2}}A_{d-1}
\xrightarrow{\alpha_{d-1}}A_d\xrightarrow{\alpha_d}A_{d+1}\xrightarrow{\alpha_{d+1}}\Sigma^d A_0.$$
where each $A_i\in\A$. It is clear that $\C$ is extension closed in $\C$.

\begin{lemma}\emph{\cite[Theorem 3.7]{l1}}\label{y1}
Let $\C$ be a $(d+2)$-angulated category and $\X\subseteq\A$ be two subcategories of $\C$.
If $(\A,\A)$ is an $\X$-mutation pair and $\A$ is extension closed,
then the quotient category $\A/\X$ is a $(d+2)$-angulated category.
\end{lemma}

\begin{lemma}\label{y2}
Let $\C$ be a $(d+2)$-angulated category with a Serre functor $\mathbb{S}$,
$\X$ be a subcategory of $\C$ and let
\begin{equation}\label{equ1}
\begin{array}{l}
\xymatrix {A_0 \xrightarrow{~\alpha_0~}A_1 \xrightarrow{~\alpha_1~} A_2 \xrightarrow{~\alpha_2~} \cdots
  \xrightarrow{~\alpha_{d - 1}~} A_d \xrightarrow{~\alpha_{d}~} A_{d+1}\xrightarrow{~\alpha_{d+1}~} \Sigma^d A_0}
\end{array}
\end{equation}
be a $(d+2)$-angle. Assume that  $\tau_d\X=\X$. Then  $\alpha_0$ is an $\X$-monic if and only if $\alpha_d$ is an $\X$-epic.
\end{lemma}

\proof  Apply the functor $\Hom(X,-)$ with $X\in\X$ to the $(d+2)$-angle (\ref{equ1}), we have the following
long exact sequence:
$$\Hom(X,A_d)\xrightarrow{~}\Hom(X,A_{d+1})\xrightarrow{~}\Hom(X,\Sigma^dA_0)\xrightarrow{~}\Hom(X,\Sigma^dA_1)$$
For $\alpha_d\colon A_d\to A_{d+1}$ to be an $\X$-epic is the same as for the first arrow in the long
exact sequence always to be epimorphism. This is the same as for the second arrow
always to be zero, which is again the same as for the third arrow always to be
injective.

Using Serre duality, the third arrow can be identified with
$$\Hom(A_0,\Sigma^{-d}\mathbb{S}X)^{\ast}\longrightarrow\Hom(A_1,\Sigma^{-d}\mathbb{S}X)^{\ast}$$
which is monomorphism if and only if
$$\xymatrix{\Hom(A_1,\Sigma^{-d}\mathbb{S}X)\ar[r]\ar@{=}[d]&\Hom(A_0,\Sigma^{-d}\mathbb{S}X)\ar@{=}[d]\\
\Hom(A_1,\tau_dX)\ar[r]&\Hom(A_0,\tau_dX)}$$
is epimorphism.
For this always to be  an epimorphism is the same as for $\alpha_0\colon A_0\to A_1$ to be a $(\tau\X)$-monic, that is, an $\X$-monic. \qed

\begin{theorem}\label{main3}
Let $\C$ be a $(d+2)$-angulated category with Serre functor $\mathbb{S}$, and $\X$ be a strongly functorially finite
subcategory of $\C$. Then the following statements are equivalent:
\begin{itemize}
\item[\rm (1)] $(\C,\C)$ is an $\X$-mutation pair.

\item[\rm (2)] The quotient category $\C/\X$ is a $(d+2)$-angulated category.

\item[\rm (3)] $\tau_d\X=\X$.
\end{itemize}
\end{theorem}

\proof (1) $\Longrightarrow$ (2). By Lemma \ref{y1}, we have that $\C/\X$ is a $(d+2)$-angulated category.

(2) $\Longrightarrow$ (1). Since $\X$ is strongly covariantly finite, for any object $A\in\C$,
by Lemma \ref{lem0}, there exists
a $(d+2)$-angle in $\C$:
$$\xymatrix {A\xrightarrow{~\alpha_0~}X_1 \xrightarrow{~\alpha_1~} X_2 \xrightarrow{~\alpha_2~} \cdots
  \xrightarrow{~\alpha_{d - 1}~} X_d \xrightarrow{~\alpha_{d}~}\mathbb{H}A\xrightarrow{~\alpha_{d+1}~} \Sigma^d A}$$
where $X_1,\cdots,X_d\in\X$ and $\alpha_0$ is a left $\X$-approximation of $A$.
Since  $\X$ is strongly contravariantly finite, for the object $\mathbb{H}A\in\C$, by Lemma \ref{duallemma}, there exists
a $(d+2)$-angle in $\C$:
$$\xymatrix {\mathbb{L}\mathbb{H}A\xrightarrow{~\beta_0~}Y_1 \xrightarrow{~\beta_1~} Y_2 \xrightarrow{~\beta_2~} \cdots
  \xrightarrow{~\beta_{d - 1}~} Y_d \xrightarrow{~\beta_{d}~}\mathbb{H}A\xrightarrow{~\beta_{d+1}~} \Sigma^d \mathbb{L}\mathbb{H}A}$$
where $Y_1,\cdots,Y_d\in\X$ and $\beta_d$ is a right $\X$-approximation of $\mathbb{H}A$.

Since $\beta_d$ is a right $\X$-approximation of $\mathbb{H}A$ and $X_d\in\X$, then
there exists a morphism $\varphi_d\colon X_d\to Y_d$ such that
$\beta_d\varphi_d=\alpha_d$.
Thus we have the following morphisms of $(d+2)$-angles
$$\xymatrix@C=1.2cm{
A\ar[r]^{\alpha_0}\ar@{-->}[d]^{\varphi}&X_1 \ar[r]^{\alpha_1}\ar@{-->}[d]^{\varphi_1}  & X_2 \ar[r]^{\alpha_2}\ar@{-->}[d]^{\varphi_2} & \cdots \ar[r]^{\alpha_{d-1}}& X_d \ar[r]^{\alpha_d\;\;}\ar[d]^{\varphi_{d}}&\mathbb{H}A\ar[r]^{\alpha_{d+1}}\ar@{=}[d]& \Sigma^dA \ar@{-->}[d]^{\Sigma^d\varphi}\\
\mathbb{L}\mathbb{H}A\ar[r]^{\beta_0}&Y_1 \ar[r]^{\beta_1} & Y_2\ar[r]^{\beta_2} & \cdots \ar[r]^{\beta_{d-1}} & Y_d \ar[r]^{\beta_d\;\;}&\mathbb{H}A\ar[r]^{\beta_{d+1}\;}& \Sigma^d\mathbb{L}\mathbb{H}A.}$$
Since $\C/\X$ is a $(d+2)$-angulated category, then the morphism
$\overline{\varphi}\colon A\to\mathbb{L}\mathbb{H}A$ is invertible in $\C/\X$.
It follows that there exists a morphism $\varphi'\colon\mathbb{L}\mathbb{H}A\to A$ such that
$\overline{\varphi'}\circ\overline{\varphi}={\rm Id}_A$ and then ${\rm Id}_A-\varphi'\varphi$ factorizes through an object in $\X$.
So it factorizes through the left$\X$-approximation $\alpha_0\colon A\to X_1$ of $A$.
Namely there exists a morphism $u\colon X_1\to A$ such that ${\rm Id}_A-\varphi'\varphi=u\alpha_0$.
Note that $({\rm Id}_A-\varphi'\varphi)\circ\Sigma^{-d}\alpha_{d+1}=u(\alpha_0\circ\Sigma^{-d}\alpha_{d+1})=0$
and then $$\Sigma^{-d}\alpha_{d+1}=(\varphi'\varphi)\circ\Sigma^{-d}\alpha_{d+1}
=\varphi'\circ\Sigma^{-d}(\beta_{d+1}).$$
Hence $\alpha_{d+1}=\Sigma^d\varphi'\circ\beta_{d+1}$
and we also have the following morphisms of $(d+2)$-angles
$$\xymatrix@C=1.2cm{
\mathbb{L}\mathbb{H}A\ar[r]^{\beta_0}\ar[d]^{\varphi'}&Y_1 \ar[r]^{\beta_1}\ar@{-->}[d]^{\phi_1}  & Y_2 \ar[r]^{\beta_2}\ar@{-->}[d]^{\phi_2} & \cdots \ar[r]^{\beta_{d-1}}& Y_d \ar[r]^{\beta_d\;\;}\ar@{-->}[d]^{\phi_{d}}&\mathbb{H}A\ar[r]^{\beta_{d+1}}\ar@{=}[d]& \Sigma^d\mathbb{L}\mathbb{H}A \ar[d]^{\Sigma^d\varphi'}\\
A\ar[r]^{\alpha_0}&X_1 \ar[r]^{\alpha_1} & X_2\ar[r]^{\alpha_2} & \cdots \ar[r]^{\alpha_{d-1}} & X_d \ar[r]^{\alpha_d\;\;}&\mathbb{H}A\ar[r]^{\alpha_{d+1}\;}& \Sigma^dA.}$$
Now we show that $\alpha_d$ is a right $\X$-approximation of $\mathbb{H}A$.
Indeed, for any morphism
$a\colon X\to \mathbb{H}A$ with $X\in\X$, since $\beta_d$ is a right $\X$-approximation of $\mathbb{H}A$,
there exists a morphism $b\colon X\to Y_d$ such that $\beta_db=a$.
It follows that $a=\beta_db=\alpha_d(\phi_db)$.
This shows that $\alpha_d$ is a right $\X$-approximation of $\mathbb{H}A$.

Similarly, one can show that for any object $C\in\C$,
there exists
a $(d+2)$-angle in $\C$:
$$\xymatrix {D\xrightarrow{~\gamma_0~}Z_1 \xrightarrow{~\gamma_1~} Z_2 \xrightarrow{~\gamma_2~} \cdots
  \xrightarrow{~\gamma_{d - 1}~} Z_d \xrightarrow{~\gamma_{d}~}C\xrightarrow{~\gamma_{d+1}~} \Sigma^d A}$$
where $Z_1,\cdots,Z_d\in\X$, $\gamma_0$ is a left $\X$-approximation and $\alpha_d$ is a right
$\X$-approximation.

This shows that $(\C,\C)$ is an $\X$-mutation pair.
\medskip

(1) $\Longrightarrow$ (3).  It suffices to show that $\X\subseteq\mathbb{S}\Sigma^{-d}\X$ and $\Sigma^d\mathbb{S}^{-1}\subseteq\X$. We only show that $\X\subseteq\mathbb{S}\Sigma^{-d}\X$, dually, we can show that $\Sigma^d\mathbb{S}^{-1}\subseteq\X$.

For any indecomposable object $X\in\X$, then $\Sigma^d\mathbb{S}^{-1}X$ is also indecomposable object. Since $\C$ has a Serre functor $\mathbb{S}$, by Theorem \ref{main2}, there exists an
Auslander-Reiten $(d+2)$-angle
$$\xymatrix {X\xrightarrow{~\alpha_0~}A_1 \xrightarrow{~\alpha_1~} \cdots
  \xrightarrow{~\alpha_{d-2}~}A_{d-1}\xrightarrow{~\alpha_{d-1}~}A_d \xrightarrow{~\alpha_d~}\Sigma^d\mathbb{S}^{-1}X\xrightarrow{~\alpha_{d+1}~} \Sigma^d X}$$
Since $(\C,\C)$ is an $\X$-mutation pair, then for the object $\Sigma^d\mathbb{S}^{-1}X$, there exists
a $(d+2)$-angle
$$\xymatrix {B\xrightarrow{~\beta_0~}X_1 \xrightarrow{~\beta_1~} X_2 \xrightarrow{~\beta_2~} \cdots
  \xrightarrow{~\beta_{d - 1}~} X_d \xrightarrow{~\beta_{d}~}\Sigma^d\mathbb{S}^{-1}X\xrightarrow{~\beta_{d+1}~} \Sigma^d B}$$
where $X_1,\cdots,X_d\in\X$, $\beta_0$ is a left $\X$-approximation and $\beta_d$ is a right $\X$-approximation.

If $\Sigma^d\mathbb{S}^{-1}X\notin\X$, then $\alpha_d\colon X_d\to\Sigma^d\mathbb{S}^{-1}X$ is not a split epimorphism.
Thus there exists a morphism $u\colon X_n\to A_n$ such that $\beta_d=\alpha_du$.
It follows that $\alpha_{d+1}\beta_d=\alpha_{d+1}\alpha_du=0$.
So there exists a morphism $v\colon \Sigma^dB\to\Sigma^dX$ such that
$v\beta_{d+1}=\alpha_{d+1}$.
Since $\beta_0$ is a left $\X$-approximation and $X\in\X$, there exists a morphism
$w\colon X_1\to X$ such that $w\beta_0=\Sigma^{-d}v$ and then $v=\Sigma^dw\circ\Sigma^d\beta_0$.
Hence $\alpha_{d+1}=v\beta_{d+1}=\Sigma^dw\circ(\Sigma^d\beta_0\circ\beta_{d+1})=0$.
This is a contradiction.

Thus we obtain $\Sigma^d\mathbb{S}^{-1}X\in\X$ and $X\in\mathbb{S}\Sigma^{-d}\X$.
Since $\C$ is a Krull-Schmidt category we infer that $\X\subseteq\mathbb{S}\Sigma^{-d}\X$.

\medskip
(3) $\Longrightarrow$ (1). It follows from Lemma \ref{y2}.  \qed
\medskip

In Theorem \ref{main3}, if $d=1$, we have the following result.

\begin{corollary}\emph{\cite[Theorem 3.3]{j}}
Let $\C$ be a triangulated category with Serre functor $\mathbb{S}$, and $\X$ be a functorially finite subcategory of $\C$. Then the quotient
category $\C /\X$ is a triangulated category if and only if $\tau\X=\X$.
\end{corollary}

Now we give an example illustrating our main result in this section.

\begin{example}
This example comes from \cite{l1}.
Let $\T=D^b(kQ)/\tau^{-1}[1]$ be the cluster category of type $A_3$,
where $Q$ is the quiver $1\xrightarrow{~\alpha~}2\xrightarrow{~\beta~}3$,
$D^b(kQ)$ is the bounded derived category of finite generated modules over $kQ$, $\tau$ is the Auslander-Reiten translation and $[1]$ is the shift functor
of $D^b(kQ)$. Then $\T$ is a $2$-Calabi-Yau triangulated category. Its shift functor is also denoted by $[1]$.

We describe the
Auslander-Reiten quiver of $\T$ in the following:
$$\xymatrix@C=0.6cm@R0.3cm{
&&P_1\ar[dr]
&&S_3[1]\ar[dr]
&&\\
&P_2 \ar@{.}[rr] \ar[dr] \ar[ur]
&&I_2 \ar@{.}[rr] \ar[dr] \ar[ur]
&&P_2[1]\ar[dr]\\
S_3\ar[ur]\ar@{.}[rr]&&S_2\ar[ur]\ar@{.}[rr]
&&S_1\ar[ur]\ar@{.}[rr]&&P_1[1]
}
$$
It is straightforward to verify that $\C:=\add(S_3\oplus P_1\oplus S_1)$ is a $2$-cluster tilting subcategory of $\T$. Moreover, $\C[2]=\C$.  By \cite[Theorem 1]{gko}, we know that $(\C,[2])$ is a $4$-angulated category with a Serre functor $\mathbb{S}=[2]$.
Let $\X=\add(S_3\oplus S_1)$.
Then the $4$-angle
$$P_1\xrightarrow{~~}S_1\xrightarrow{~~}S_3\xrightarrow{~~}P_1\xrightarrow{~~}P_1[2]$$
shows that $\X$ is strongly functorially finite subcategory of $\C$ and $(\C,\C)$ is an $\X$-mutation pair.
By Theorem \ref{main3}, we obtain that $\C/\X$ is 4-angulated category and $\tau_2\X=\X$.
\end{example}

%\section*{Acknowledgments}
%
%The author is highly grateful to the referee for his/her valuable comments and suggestions which led to improvements of this paper.

\textbf{Panyue Zhou}\\
College of Mathematics, Hunan Institute of Science and Technology, 414006 Yueyang, Hunan, P. R. China.\\
E-mail: \textsf{panyuezhou@163.com}


\begin{thebibliography}{99}
\addtolength{\itemsep}{-0.18em}

\bibitem[AHBT]{ahbt} E. Arentz-Hansen, P.  Bergh, M. Thaule. The morphism axiom for $N$-angulated categories. Theory Appl. Categ. 31(18):  477--483, 2016.

\bibitem[AR]{ar} M. Auslander, I. Reiten. Representation theory of Artin algebras. IV. Invariants given by almost split sequences. Comm. Algebra 5(5): 443--518, 1977.

\bibitem[B]{b}
R. Bocklandt. Graded Calabi--Yau algebras of dimension $3$, with an appendix ``The
signs of Serre duality" by M. Van den Bergh. J. Pure Appl. Algebra 212(1):
14--32, 2008.

\bibitem[BK]{bk} A. Bondal, M. Kapranov. Representable functors, Serre functors, and mutations.
Mathematics of the USSR-Izvestiya, 35(3): 519--541, 1990.

\bibitem[BT1]{bt1} P. Bergh, M. Thaule. The axioms for $n$-angulated categories. Algebr. Geom. Topol. 13(4): 2405--2428, 2013.

\bibitem[BT2]{bt2} P. Bergh, M. Thaule.
The Grothendieck group of an $n$-angulated category.
J. Pure Appl. Algebra, 218(2): 354--366, 2014.



\bibitem[F]{f}
F. Fedele. Auslander-Reiten $(d+2)$-angles in subcategories and a $(d+2)$-angulated generalisation of a theorem by Br\"{u}ning. J. Pure Appl. Algebra, 223(8): 3554--3580, 2019.

\bibitem[GKO]{gko} C. Geiss, B. Keller, S. Oppermann. $n$-angulated categories.  J. Reine Angew. Math.  675: 101--120, 2013.

\bibitem[IY]{iy} O. Iyama, Y. Yoshino. Mutations in triangulated categories and rigid Cohen-Macaulay modules. Invent. Math. 172(1): 117--168, 2008.

\bibitem[Ja]{ja} G. Jasso. $n$-abelian and $n$-exact categories. Math. Z. 283(3-4):  703--759, 2016.


\bibitem[J]{j} P. J{\o}rgensen.
Quotients of cluster categories.
Proc. Roy. Soc. Edinburgh Sect. A 140(1):  65--81, 2010.

\bibitem[L1]{l1} Z. Lin. $n$-angulated quotient categories induced by mutation pairs.
Czechoslovak Math. J. 65(4): 953--968, 2015.

\bibitem[L2]{l2}  Z. Lin. Right $n$-angulated categories arising from covariantly finite
subcategories. Comm. Algebra, 45(2): 828--840, 2017.

\bibitem[OT]{ot} S. Oppermann, H. Thomas.  Higher-dimensional cluster combinatorics and representation theory. J. Eur. Math. Soc. 14(6): 1679--1737, 2012.


\bibitem[RVB]{rv} I. Reiten, M. Van den Bergh.
Noetherian hereditary abelian categories satisfying Serre duality.
J. Amer. Math. Soc. 15(2): 295--366, 2002.

\bibitem[Z]{z}
P. Zhang. Triangulated categories and derived categories. Chinese. Science Press, 2015.
\end{thebibliography}
\end{document}